\PassOptionsToPackage{hyphens}{url}
\PassOptionsToPackage{hypcap=false}{caption}
\documentclass[a4paper,cleveref,autoref,thm-restate]{article}

\usepackage{amsmath}
\usepackage{amsfonts}
\usepackage{amssymb}
\usepackage{graphicx}
\usepackage{enumerate}
\usepackage{amsthm}
\usepackage[a4paper,includeheadfoot,margin=1in]{geometry}
\usepackage{listings}
\usepackage{tikz}
\usepackage{hyperref}
\usepackage{xcolor}
\usepackage{mathtools}
\usepackage{pdfpages}
\usepackage{multirow}
\usepackage{hhline}
\usepackage{caption}
\usepackage{subcaption}
\usepackage{colortbl}
\usepackage{ragged2e}
\usepackage{parskip}
\usepackage{tabularx}
\usepackage{overpic}
\usepackage{wrapfig}
\usepackage{tcolorbox}
\usepackage{algorithm2e}
\usepackage{seqsplit}

\newcommand{\regina}{\textit{Regina} }

\newcommand{\is}[1]{{\ttfamily\seqsplit{#1}}}

\newcommand{\pt}{ \ensuremath{\{ \mathrm{pt.} \}} }

\newcommand{\svert}{ \ensuremath{\,\vert\,} }

\DeclareMathOperator\lk{lk}

\DeclareTextFontCommand{\bemph}{\boldmath\bfseries}

\newtheorem{theorem}{Theorem}

\newtheorem{conjecture}{Conjecture}

\newtheorem{remark}[theorem]{Remark}

\newtheorem{definition}[theorem]{Definition}

\title{Practical Software for Triangulating and Simplifying 4-Manifolds}

\author{Rhuaidi Antonio Burke$^1$}
\date{
{\itshape\small{$^1$School of Mathematics and Physics, University of Queensland, Brisbane QLD 4072, Australia}}\\
{\ttfamily\small{rhuaidi.burke@uq.edu.au}}\\[2ex]%
\today
}

\begin{document}

\maketitle

\begin{abstract}\noindent
Dimension 4 is the first dimension in which exotic smooth manifold pairs appear -- manifolds which are topologically the same but for which there is no smooth deformation of one into the other. Whilst smooth and triangulated 4-manifolds do coincide, comparatively little work has been done towards gaining an understanding of smooth 4-manifolds from the discrete and algorithmic perspective. In this paper we introduce new software tools to make this possible, including a software implementation of an algorithm which enables us to build triangulations of 4-manifolds from Kirby diagrams, as well as a new heuristic for simplifying 4-manifold triangulations. Using these tools, we present new triangulations of several bounded exotic pairs, corks and plugs (objects responsible for ``exoticity''), as well as the smallest known triangulation of the fundamental $K3$ surface. The small size of these triangulations benefit us by revealing fine structural features in 4-manifold triangulations.
\end{abstract}


\paragraph{Acknowledgements}We thank the referees for their helpful comments.
The author was supported by an Australian Government Research Training Program Scholarship.

\section{Introduction}
\label{sec:intro}
In dimensions $\leq 3$, every topological manifold admits a unique smooth structure. In dimensions $\geq 4$ however, this is no longer the case --- there exist manifolds which are \bemph{homeomorphic} (they represent the same topological manifold), but not \bemph{diffeomorphic} (they represent distinct smooth manifolds). Such manifolds are called \bemph{exotic}. 

One of the great remaining open problems of classical topology is the smooth 4-dimensional Poincar\'e conjecture (SPC4), which asserts that all 4-manifolds homeomorphic to the 4-sphere are diffeomorphic. In other words, SPC4 asks whether or not there exist exotic 4-spheres.

By the work of Cairns \cite{MR1563139, MR0149491} and Whitehead \cite{MR0002545}, every smooth manifold can be uniquely triangulated (a smooth structure uniquely determines a piecewise-linear (PL) structure). In dimensions $\leq 6$, the converse also holds, so PL manifolds admit a unique smooth structure \cite{MR0415630,MR0121804}.  
In particular then, results which hold for smooth 4-manifolds also hold for PL 4-manifolds, and so we may move between the two settings at our discretion.

As such, one might hope to gain an insight into smooth 4-manifolds by studying their PL structures, specifically in terms of triangulations, which are better suited to computational techniques.
Despite this, comparatively little work has been done from this perspective. 

One reason for this may stem from the fact that in contrast to dimension three, where algorithms to solve many problems are known (though in many cases are intractable), in dimension four many of the same problems become undecidable. For example, in dimension three there exists an algorithm to decide whether two arbitrary triangulations represent the same topological manifold (despite the best known algorithm having complexity bounded by a tower of exponentials) \cite{Kuperberg-M3Homeo,Lackenby}. In dimension four however, any finitely presented group can be the fundamental group of a 4-manifold, from which the problem of deciding whether two arbitrary 4-manifolds are homeomorphic (let alone diffeomorphic) becomes equivalent to solving the word problem on finitely presented groups, which is undecidable \cite{MR0081851,MR0097793,MR0110743}. 

As a consequence, we can often only hope for heuristics which, for as many cases as possible, give the correct answer, in as short a time as possible. The existence of such heuristics illustrate the difference between what can be shown in theory versus what is possible in practice, and which motivates the contents of this paper. Of course, there are theoretical restrictions on the effectiveness of any heuristic (for example there is no computable upper bound on the number of local moves needed to simplify to a minimal triangulation, since this would lead to a violation of the aforementioned word problem).

\bemph{Original Contributions} Regarding the proposed program of understanding smooth, and in particular, exotic 4-manifolds from a computational perspective, and motivated by SPC4, we should ask: \textit{what tools and techniques will be needed for such an undertaking?} The most obvious starting point would be to build up a catalogue of examples to analyse, consisting of triangulations of different, simple, exotic pairs. In particular, if we are interested in eventually attacking SPC4, we desire examples which are closed, orientable, simply-connected, and in some sense ``small'' (either topologically, for example in terms of the Euler characteristic, or in terms of the number of simplices in a triangulation, or indeed both). 

Concerning a catalogue of examples, until now there has been just one readily available pair of ``exotic'' triangulations, realising a pair due to Kreck \cite{kreck} and triangulated by Benedetti and Lutz \cite{BenedettiLutz}. The manifolds of the pair are formed from connect sums of standard 4-manifolds, and are non-orientable and not simply-connected. Moreover, the triangulations are presented as simplicial complexes and so are incredibly large (having 460 and 518 4-simplices respectively), in principle making a detailed analysis of their combinatorics difficult due to a lack of effective simplification techniques, until now.

Whilst one typically associates the term ``triangulation'' with ``simplicial complex'', we work here instead with \bemph{generalised triangulations} (specifically, unordered $\Delta$-complexes). These allow for far more efficient triangulations compared to simplicial complexes since, for example, faces of the same simplex can be identified together.

We typically want as small a triangulation as possible in order, for example, to recognise important or otherwise interesting structures within the triangulation. As such, we require an effective means of simplifying triangulations. Existing simplification heuristics tend to be predominantly tailored for either  3-manifold triangulations, or triangulations presented as simplicial complexes (and which retain the simplicial structure), whilst simplifying 4-manifold triangulations remains a challenge. 

In this article, we address the above concerns, presenting: 
\begin{itemize}
	\item In Section \ref{sec:dgt}: Software which allows for the first time, the ability to produce triangulations of 4-manifolds from a commonly used smooth description (Kirby diagrams).
	\item In Section \ref{sec:uds}: A new and effective heuristic for simplifying 4-manifold triangulations.
	\item In Section \ref{sec:results}, using these software tools we obtain:
	\begin{itemize}
		\item The first examples of small triangulations (ranging from 10 to 26 simplices) of several bounded, simply-connected, orientable, exotic pairs.
		\item Small triangulations of objects which are directly responsible for the ``exoticity'' in certain exotic pairs (corks and plugs).
		\item The current smallest known triangulation of the $K3$ surface (one of four ``fundamental'' simply-connected 4-manifolds).
	\end{itemize}
\end{itemize}
As a consequence of the above, we are for the first time in a position to begin effectively analysing the combinatorics of both standard and exotic 4-manifold triangulations, which we illustrate in Section \ref{sec:results}. Finally, in Section \ref{sec:ongoing}, we discuss ongoing work to produce our first triangulations of \emph{closed}, simply-connected, orientable, exotic 4-manifolds.

\section{Preliminaries}
\label{sec:bgtheory}

\subsection{Links and 4-Manifolds}
\label{sec:kirby-diagrams}
We assume a basic familiarity with knot theory (writhes, etc.); for details see \cite{Adams-KnotBook}. Henceforth all manifolds are assumed to be orientable. In the smooth setting, we primarily work with 4-manifolds via \bemph{handle} decompositions \cite{Akb-4mflds,MR1707327}.

\begin{definition}
For $0\leq k\leq 4$, a \bemph{4-dimensional $k$-handle} $h^k$ is a copy of $D^k\times D^{4-k}$, attached to a smooth $4$-manifold $W$ via an embedding of the form $\varphi:\partial D^k\times D^{4-k}\to\partial W$.	
\end{definition}

\begin{figure}[h]
\centering
\includegraphics[width=\textwidth]{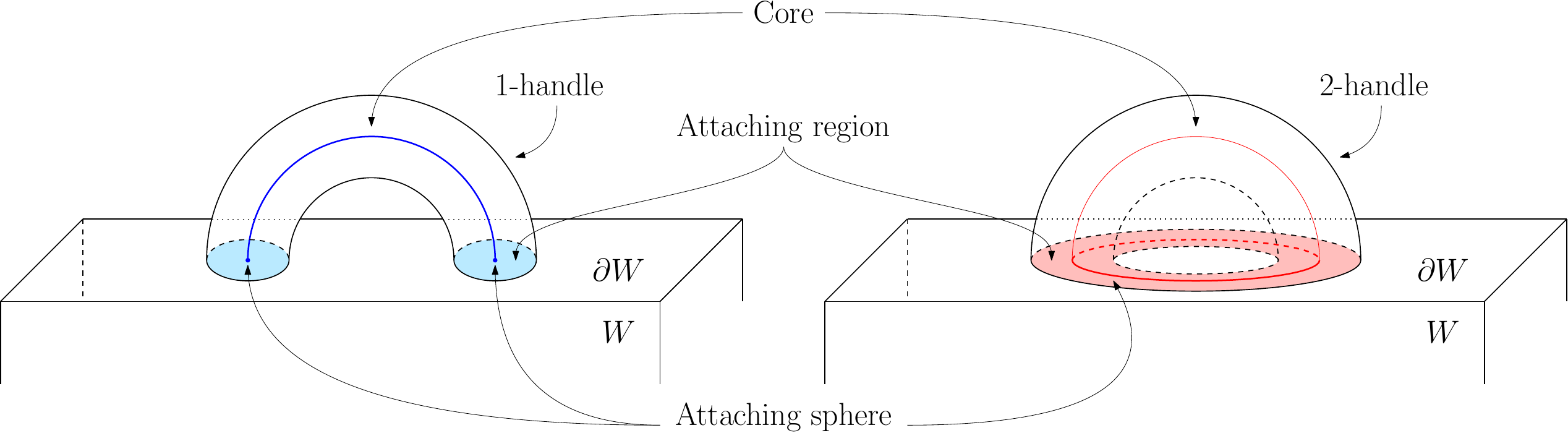}
\caption{Cartoons illustrating adding a 1-handle versus adding a 2-handle.}
\label{fig:handleCartoon}
\end{figure}

All handles are topologically just a ball, so what distinguishes a handle is \textit{how} it is attached. For example, attaching a 1-handle can be thought of as a ``rod'' with its ends attached to the base $W$, whereas a 2-handle could be thought of as a ``plate'' with its circular boundary attached to $W$ (and one imagines ``thickening'' these rods and plates up to be the appropriate dimension so that we always get a manifold). This idea is shown in Figure \ref{fig:handleCartoon}.

We typically build a 4-manifold by starting with $W=B^4$ (a 0-handle), and then attaching 1- and 2-handles. By a theorem of Laudenbach and Poenaru, 3- and 4-handles attach uniquely \cite{MR0317343} up to diffeomorphism, and so it suffices to understand how the 1- and 2-handles attach.

Given a handle decomposition of $M$, we visualise $M$ by placing ourselves in the boundary of the 0-handle, $\mathbb{S}^3\cong\mathbb{R}^3\cup\{\infty\}$, and drawing the attaching regions of the 1- and 2-handles.

\begin{figure}[h]
	\centering
	\includegraphics[width=\textwidth]{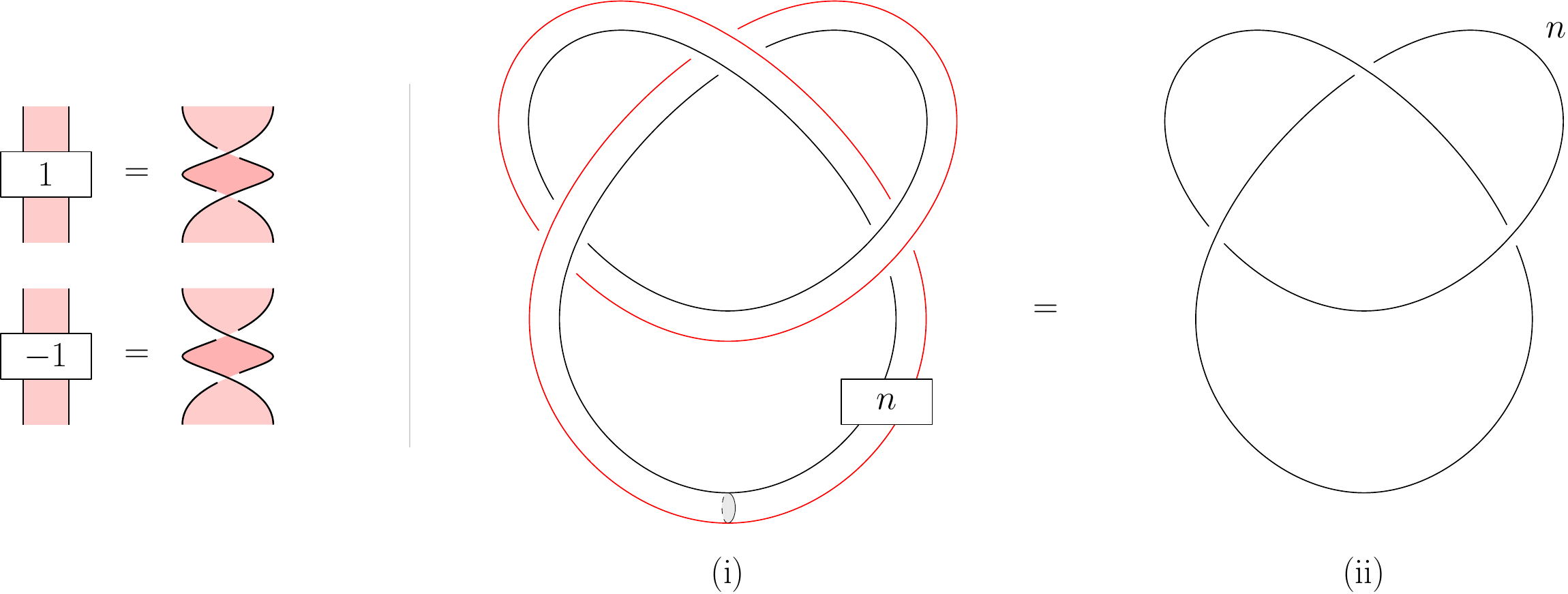}
	\caption{A framed knot.}
	\label{fig:framedKnot}
\end{figure}

Let us first consider the case of a single 2-handle attached to the 0-handle. The attaching map is of the form $\varphi:S^1\times D^2\to S^3$. Such a map is determined up to isotopy by:
\begin{enumerate}
	\item an embedding $\varphi\rvert_{S^1\times\{0\}}:S^1\times\{0\}\to S^3$, i.e.\ a \bemph{knot} $K$, and
	\item a \bemph{framing} of $\varphi(S^1)$, i.e.\ a choice of normal vector field on $K$.
\end{enumerate}
Framings are in bijection with the integers \cite{Akb-4mflds,MR1707327}. This integer is the number of times the $D^2$ factor ``twists'' around the knot. Diagrammatically this can be represented as a ``ribbon'' formed from $K$ with $n$ twists in it (Figure \ref{fig:framedKnot}(i)). More commonly however we will simply draw $K$ decorated with the integer $n$ as shown in Figure \ref{fig:framedKnot}(ii).  

\begin{figure}[h]
	\centering
	\includegraphics[width=0.35\textwidth]{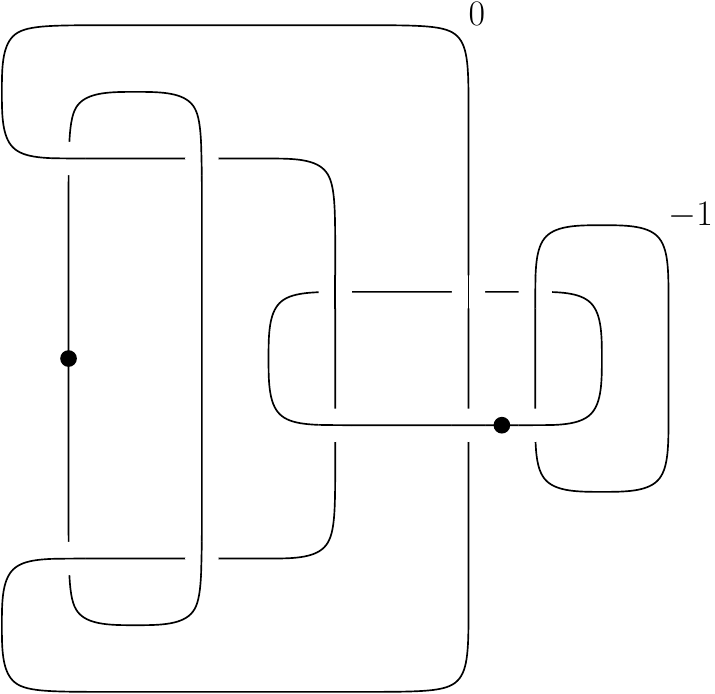}
	\caption{An example of a Kirby diagram.}
	\label{fig:kirbyDiagramExample}
\end{figure}

A 1-handle is represented diagrammatically by an unknot with a dot on it (cf.\ Figure \ref{fig:kirbyDiagramExample}). Arcs of 2-handles passing through a dotted circle are understood to be running ``over'' the attaching region of the 1-handle. See Appendix \ref{app:1handleNotation} for a technical explanation of this notation.

With the above in mind, we can visualise a smooth 4-manifold by a decorated, $\ell$-component link diagram $L$, comprised of (i) $d$ dotted, mutually unlinked unknots, representing the 1-handles; and (ii) $\ell-d$ knots (possibly linking each other and/or the dotted components), with an integer $c_i$ attached to each component $L_i$ ($d+1\leq i\leq \ell$), representing the 2-handles and their associated framings. One should also specify how many, if any, 3-handles there are. This data is enough to reconstruct the 4-manifold up to diffeomorphism \cite{Akb-4mflds,MR1707327}. Such a diagram is called a \bemph{Kirby diagram}. Figure \ref{fig:kirbyDiagramExample} is an example of a typical Kirby diagram.

\subsection{Gems \& Crystallisations}
\label{sec:crystallisation-theory}
In what follows the term \bemph{graph} is used to refer to a finite multigraph without loops.

\begin{definition}
	An \bemph{$(n+1)$-coloured graph} is a pair $(\Gamma,\gamma)$, where $\Gamma=(V(\Gamma),E(\Gamma))$ is an $(n+1)$-regular graph, and $\gamma:E(\Gamma)\to\Delta_n=\{0,\ldots,n\}$ is a map which is surjective on adjacent edges (a \bemph{colouring} of the edges so that no two adjacent edges have the same colour).
\end{definition}

Let $S\subseteq\Delta_n$ and let $\Gamma_{S}$ be the graph obtained from $\Gamma$ by deleting all the edges that are not coloured by elements of  $S$. The connected components of $\Gamma_{S}$ are called \bemph{$S$-residues} of $\Gamma$.

\begin{figure}[h]
	\centering
	\includegraphics[width=0.75\textwidth]{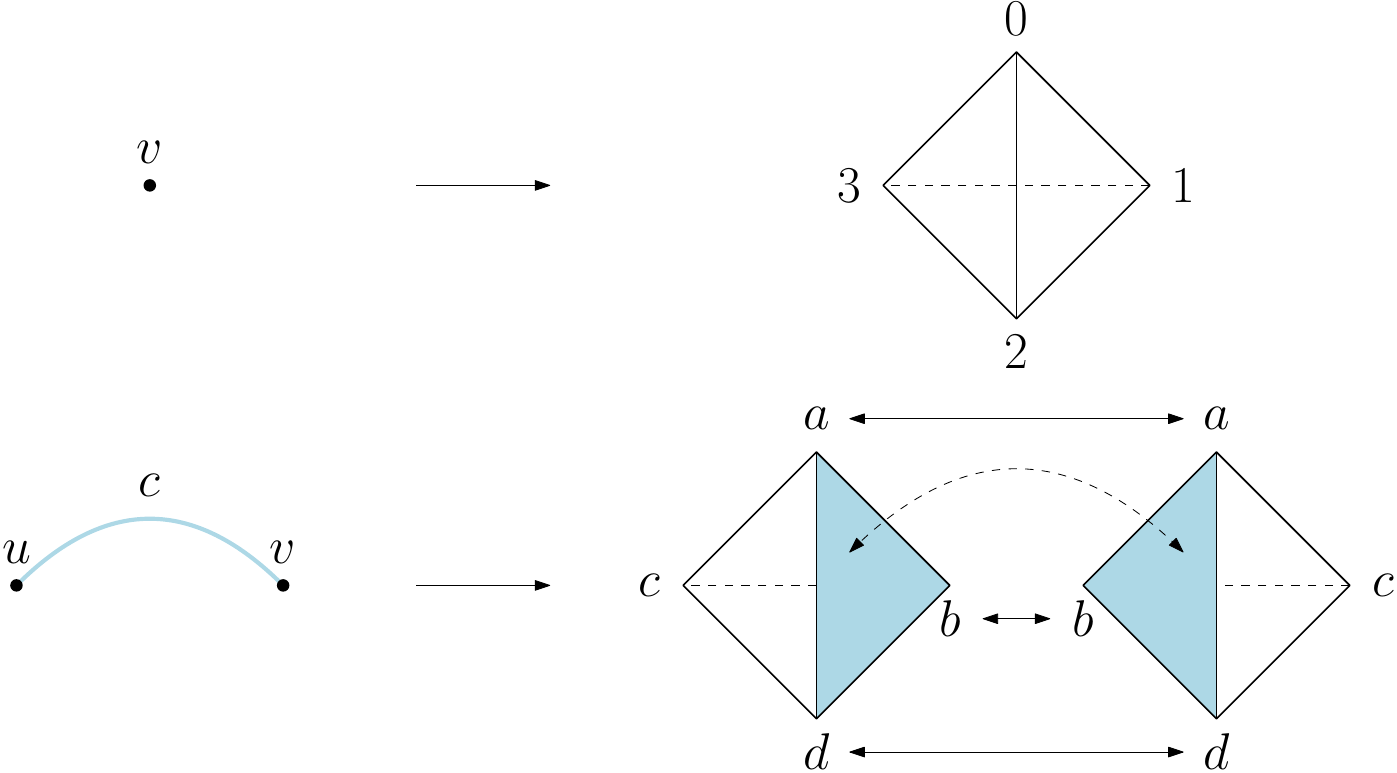}
	\caption{Constructing a triangulation from a coloured graph.}
	\label{fig:graph2tri-3dDemo}
\end{figure}

An $(n+1)$-coloured graph $(\Gamma,\gamma)$ can be used to construct an $n$-dimensional pseudo-complex $K(\Gamma)$ via the following procedure \cite{MR0867510}:
\begin{enumerate}
	\item For each $v\in V$, take an $n$-simplex $\sigma(v)$, and label its vertices with the elements of $\Delta_n$.
	\item For each $c$-coloured edge between $u$ and $v$ ($u,v\in V$), identify the $(n-1)$-faces of $\sigma(u)$ and $\sigma(v)$ which are opposite their respective $c$ labelled vertex, such that equally labelled vertices of $\sigma(u)$ and $\sigma(v)$ are identified.
\end{enumerate}
This construction is illustrated for the three-dimensional case in Figure \ref{fig:graph2tri-3dDemo}.

If $K(\Gamma)$ triangulates a PL $n$-manifold $M$, then $(\Gamma,\gamma)$ is called a \bemph{gem} (\bemph{g}raph \bemph{e}ncoded \bemph{m}anifold) representing $M$. In the case that $\partial M\neq\emptyset$, $K(\Gamma)$ triangulates the associated \emph{singular} manifold $\widehat{M}$ obtained from $M$ by taking a cone over each boundary component. Closed $n$-manifolds are a subset of singular $n$-manifolds, with $M=\widehat{M}$. For simplicity we will make no further distinction between these cases. By construction the graph $\Gamma$ represents the dual 1-skeleton of the complex $K(\Gamma)$. An $(n+1)$-coloured gem representing $M$ is called a \bemph{crystallisation} if the associated triangulation of $M$ has exactly $n+1$ vertices. 

\begin{theorem}[\cite{CasaliCristoforiGrasselli2018}]
Every compact $n$-manifold admits an $(n+1)$-coloured graph representing it.
\end{theorem}

\section{Diagrams to Graphs and Triangulations}
\label{sec:dgt}
In 2000, Casali described an algorithmic way to obtain a gem of a 4-manifold $M$ from a framed link $M=M(L,c)$ (i.e.\ a Kirby diagram that consists only of 2-handles) \cite{Casali2000}. More recently, Casali and Cristofori have extended this algorithm to the case where $L$ may now contain dotted components (i.e.\ 1-handles) \cite{CasaliCristofori2023Final}. However, until now there has been no readily available software implementation of these algorithms. We introduce here an implementation of both algorithms in a software utility called \textit{DGT} (\bemph{D}iagrams to \bemph{G}raphs \& \bemph{T}riangulations)\footnote{The author is open to suggestions for a better, snappier, name.}, as a part of the \regina software package \cite{regina}.

\begin{figure}[h]
	\centering
	\includegraphics[width=0.9\textwidth]{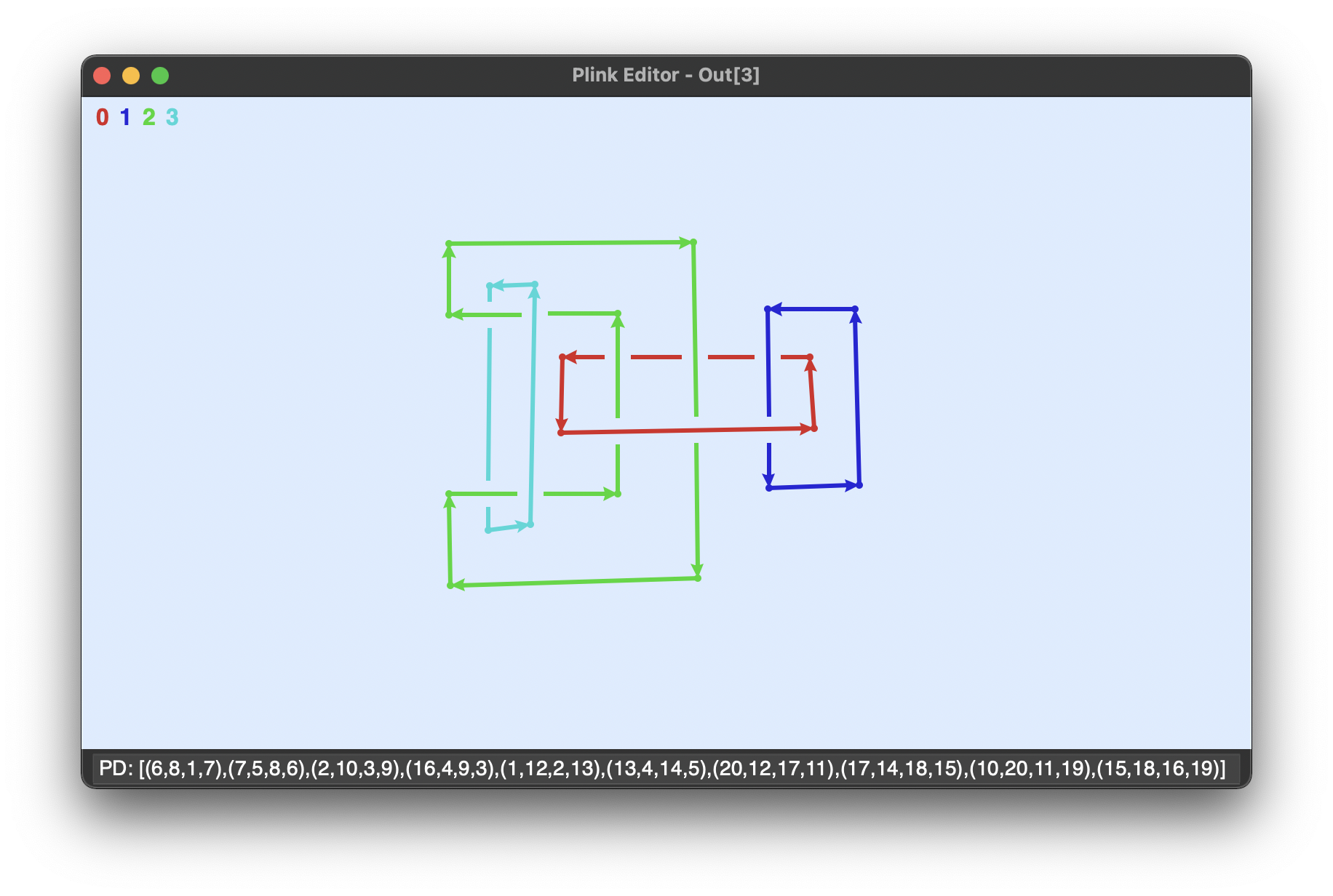}
	\caption{Snappy's PLink Editor}
	\label{fig:plinkEditor}
\end{figure}

The user provides a combinatorial encoding of the underlying link of the Kirby diagram in the form of a \bemph{Planar Diagram Code} (PD Code) (this is most easily generated by drawing the link in the \textit{PLink Editor} included as part of \textit{Snappy} \cite{SnapPy}, (see Figure \ref{fig:plinkEditor}), and a list specifying the integer framing on each component or the presence of a 1-handle. DGT then builds a 5-coloured graph realising a gem of the associated 4-manifold, from which it in turn also produces a \regina formatted triangulation.

\subsection{Construction Details}
In this section we summarise the algorithms of \cite{CasaliCristofori2023Final}. We note that our version of the algorithm in the case involving 1-handles is slightly different from what is presented in \cite{CasaliCristofori2023Final}, and is more representative of our particular implementation. 

\begin{figure}[h]
	\centering
	\includegraphics[width=0.95\textwidth]{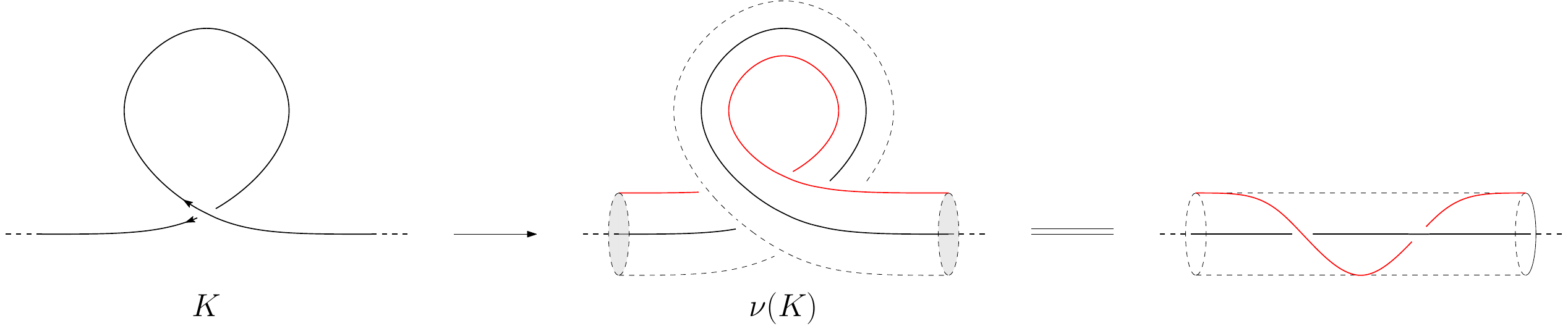}
	\caption{Adding a curl adds a twist in the normal disk bundle (i.e.\ changes the framing by $\pm 1$).}
	\label{fig:r12framing}
\end{figure}

To begin, we encode the framing of a 2-handle via the writhe of the knot (since the writhe of a knot coincides with the so-called ``blackboard framing'' of that knot). If the specified framing and the writhe do not coincide, we can add additional ``curls'' (Reidemeister 1 moves) of the appropriate sign to the knot, which change the writhe by $\pm 1$, until the writhe equals the specified framing. Figure \ref{fig:r12framing} shows how the addition of these curls corresponds to a tubular neighbourhood of the knot being twisted around the core, and hence encoding the framing as desired. DGT automatically performs this ``self-framing procedure''.

\begin{figure}[h]
	\centering
	\includegraphics[width=0.95\textwidth]{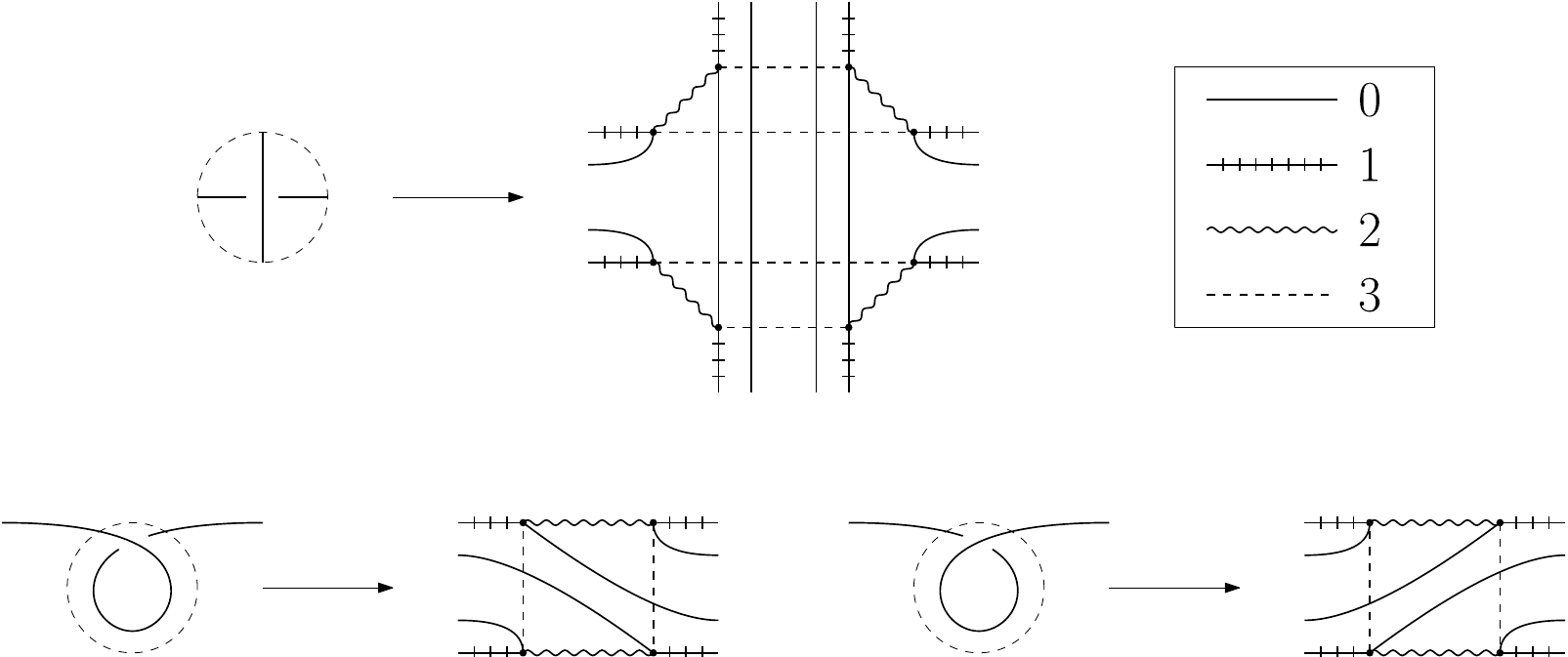}
	\caption{4-Coloured subgraphs corresponding to true crossings and (positive and negative) curls.}
	\label{fig:dgt-graphs1}
\end{figure}

To construct the 5-coloured gem $\Lambda$ representing $M$ from its Kirby diagram, the first step is to construct a 4-coloured gem $\Gamma$ representing $\partial M$\footnote{In the case that $M$ is closed, we remove a small ball and let $\partial M$ be the resulting boundary 3-sphere.}. $\Gamma$ is constructed as follows.
\begin{enumerate}
	\item For each crossing and curl in the diagram $L$, construct 4-coloured graphs as per Figure \ref{fig:dgt-graphs1}.
	\item Identify the ``hanging'' edges of each subgraph together in the ``natural'' way according to the link diagram (see \cite{Casali2000} for the precise characterisation of ``natural'').
\end{enumerate}

\begin{figure}[h]
	\centering
	\includegraphics[width=0.25\textwidth]{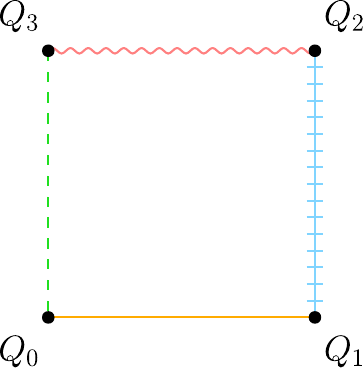}
	\caption{A quadricolour.}
	\label{fig:quadricolour}
\end{figure}

By possibly introducing a pair of cancelling curls (i.e.\ of opposite signs) into each component $L_i$ of $L$, the graph $\Gamma$ is guaranteed to contain, for each $L_i$, a copy of the subgraph shown in Figure \ref{fig:quadricolour}, called a \bemph{quadricolour}.

\begin{figure}[h]
	\centering
	\includegraphics[width=0.85\textwidth]{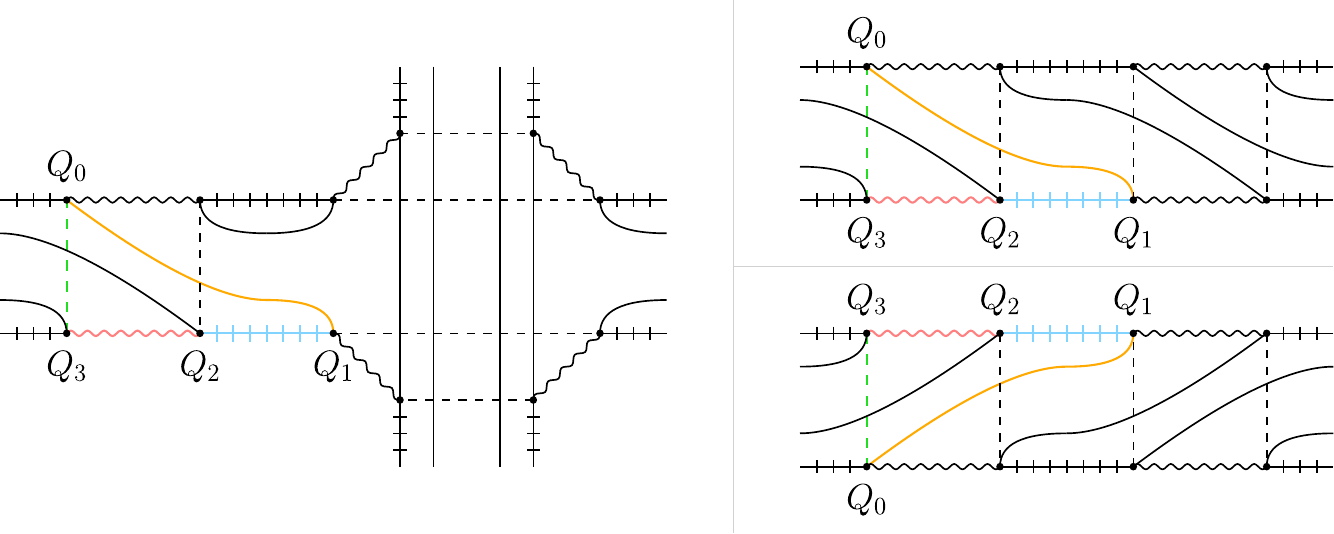}
	\caption{Guaranteed quadricolour location conditions.}
	\label{fig:quadriAppearances}
\end{figure}

Figure \ref{fig:quadriAppearances} shows that a quadricolour will appear whenever a curl is adjacent to an undercrossing or another curl of the same sign.

\begin{figure}[h]
	\centering
	\includegraphics[width=0.8\textwidth]{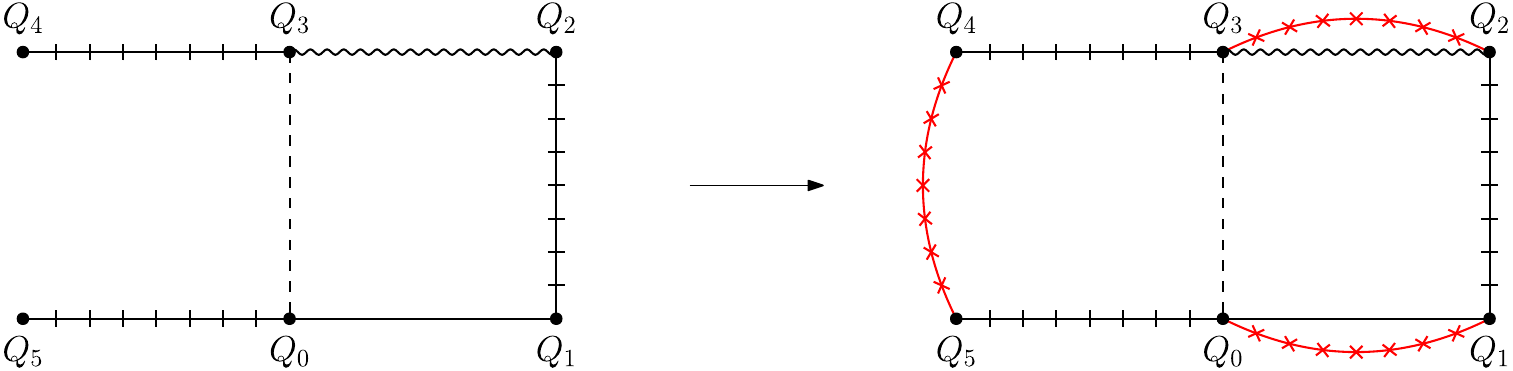}
	\caption{The quadricolour substitution operation realising 2-handle attachment.}
	\label{fig:quadriSub}
\end{figure}

For a 4-manifold $M=M(L,c)$ with no 1-handles, the procedure to build a 5-coloured graph $\Lambda$ representing $M$ can now be summarised as follows:
\begin{enumerate}
	\item Construct the 4-coloured graph $\Gamma$ representing $\partial M$.
	\item For each $L_i$, locate a quadricolour associated to that component in $\Gamma$, and add colour-4 edges according to Figure \ref{fig:quadriSub}.
	\item Add colour-4 edges to the remaining vertices of $\Gamma$ sharing a colour-1 edge.    
\end{enumerate}

In short, the operation of Figure \ref{fig:quadriSub} realises attachment of a 2-handle, whilst the doubling of 1-coloured edges realises a cone over the remaining boundary complex (see e.g.\ \cite{Casali2000}).

In theory, to locate a quadricolour for each component, one could build the 1-residue of $\Gamma$ (the connected components of which correspond to the components of $L$), and for each component, perform a search for a cyclic subgraph corresponding to a quadricolour. However we can far more efficiently locate a quadricolour as follows. 

Since a quadricolour is guaranteed to appear at either two adjacent curls of the same sign, or at a curl and the under strand of a ``true'' crossing (a non-curl crossing), we can determine where quadricolours will appear by identifying these structures directly from the diagram. Then by implementing our graphs via adjacency lists, locating a quadricolour in the graph can be made linear time using a single pass over the vertices of the graph.

\begin{figure}[h]
\centering
\includegraphics[width=0.85\textwidth]{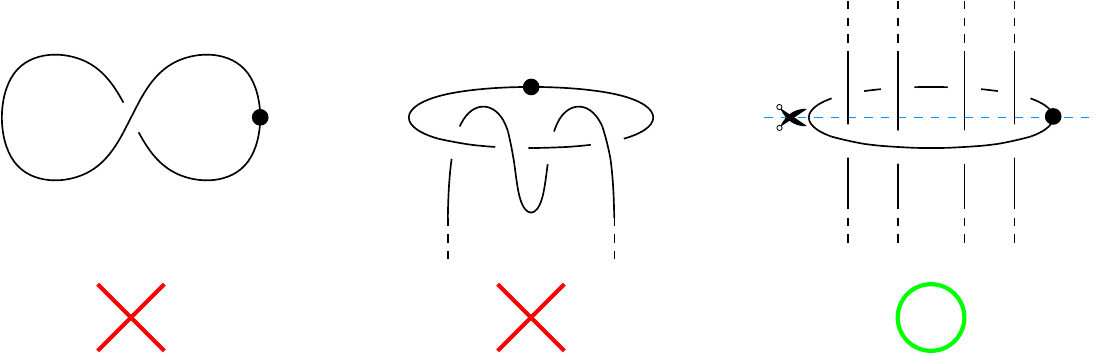}
\caption{Requirements for the handle diagram.}
\label{fig:dgt-diagramHandleReqs}
\end{figure}

Suppose we now have $d>0$ dotted components. We first require the link diagram to satisfy two conditions: (i) the 1-handles are drawn in a ``standard position'' --- that is, as unknots which visibly bound a disk --- and (ii) any arcs of framed components passing through a dotted component do so in a way such that one could ``cut'' the dotted component into two halves, one half containing only over crossings of the 2-handles, and one containing only undercrossings of the 2-handles (see Figure \ref{fig:dgt-diagramHandleReqs}). Kirby diagrams can always be isotoped to satisfy these conditions \cite{MR1707327} (though in practice this is often easier said than done).

\begin{figure}[h]
	\centering
	\includegraphics[width=0.75\textwidth]{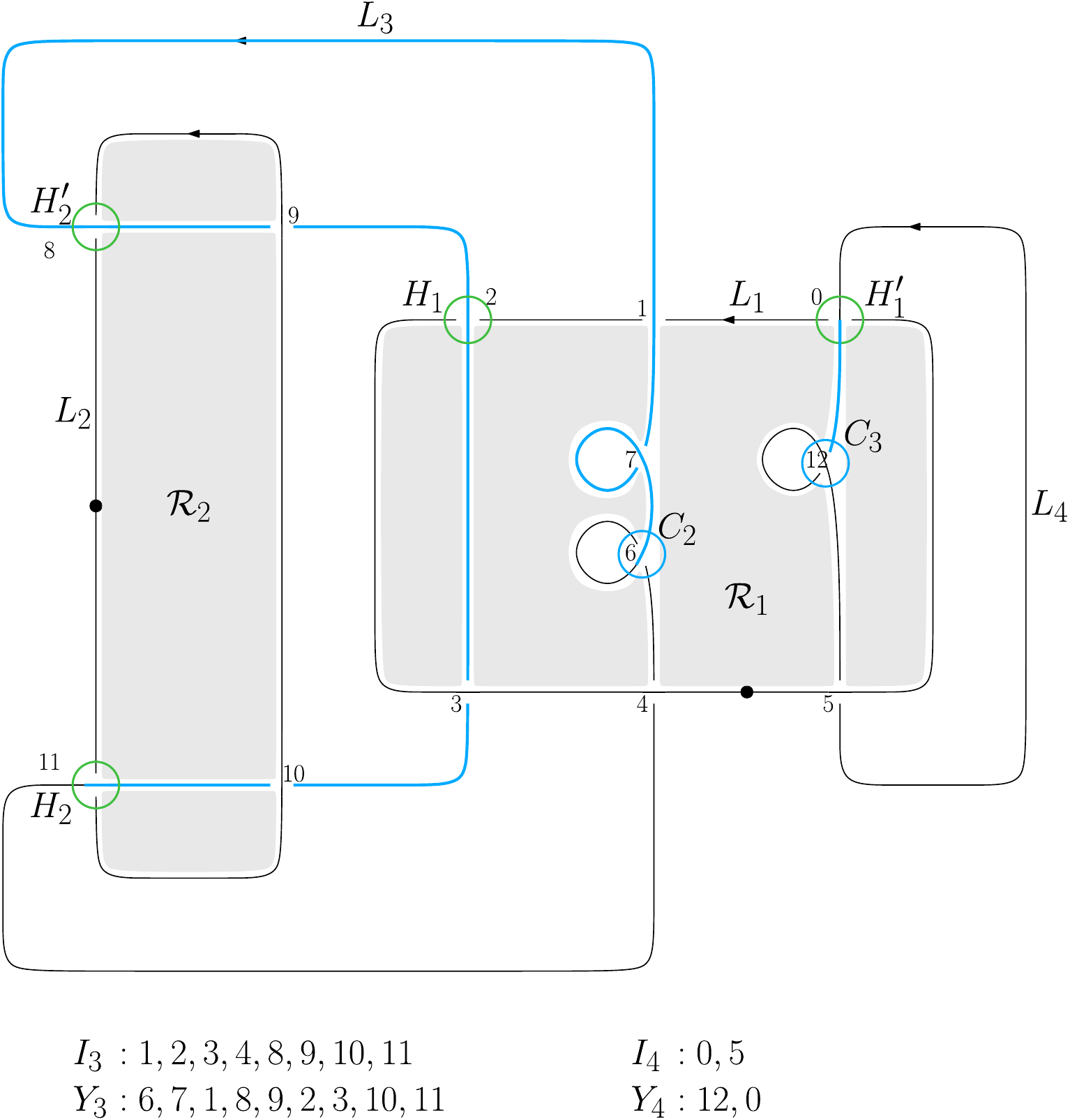}
	\caption{Example illustrating the various link data in the 1-handle construction.}
	\label{fig:dgt-1h-eg-1}
\end{figure}

Let $L_i$ denote the $i$th dotted component for $1\leq i\leq d<\ell$, and $L_j$ denote the $j$th framed component for $d+1\leq j\leq \ell$. Let $Q_j$ be the chosen quadricolour in $\Gamma$, corresponding to $L_j$. Since the location of $Q_j$ can be determined from a pair of crossings in $L$, let $(C_j,X_j)$ denote the curl and undercrossing pair in $L$ which determines $Q_j$. If $Q_j$ is a quadricolour constructed from two curls, let $X_j$ be the `other' curl in the pair. Record the following data:
\begin{itemize}
	\item For each $i\in\{1,\ldots,d\}$, let $\{H_i,H_i'\}$ be the two ``outermost'' undercrossings of $L_i$ (see Figure \ref{fig:dgt-1h-eg-1}). 
\item For each $j\in\{d+1,\ldots,\ell\}$: 
\begin{itemize}
	\item Let $I_j:=L_j\cap\bigcup_{i=1}^{d} L_i$ be the set of crossings common to both $L_j$ and $L_i$ for a fixed $j$.
	\item Starting from $C_j$, walk along $L_j$ in the direction opposite to $X_j$. Let $Y_j$ be the set of crossings encountered. Walk along $L_j$ until $Y_j$ contains all crossings in $I_{j}$ except possibly for $X_j$. Let $Y=\bigcup_j Y_j$.
\end{itemize}
\end{itemize}

In essence, this data ``finds'' the disks which the dotted components would bound in the absence of any 2-handles, thereby effectively locating the 1-handles in the diagram. This procedure is shown in Figure \ref{fig:dgt-1h-eg-1}, applied to the Kirby diagram depicted in Figure \ref{fig:kirbyDiagramExample}.

Once the 4-colour graph $\Gamma$ has been constructed, and using the data from above, record the following $\Gamma$-related data:
\begin{itemize}
	\item For each $i\in\{1,\ldots,d\}$ let $\mathcal{R}_i$ be region bounded by $L_i$. For each $i$, let $e_i$ (respectively $e_i'$) be the 1-coloured edges of $\Gamma$ parallel to the part of the arc of $L_i$ containing $H_i$ (resp.\ $H_i'$) on the side of the regions of $L$ ``merging'' into $\mathcal{R}_i$, and let $v_i$ (resp.\ $v_i'$) be its endpoint belonging to the subgraph corresponding to an undercrossing of $L_i$ (see Figure \ref{fig:dgt-1h-eg-2}).
\end{itemize}
See \cite{CasaliCristofori2023Final} for the crystallisation-theoretic details of why this data is needed.

\begin{figure}[h]
\centering
\includegraphics[width=\textwidth]{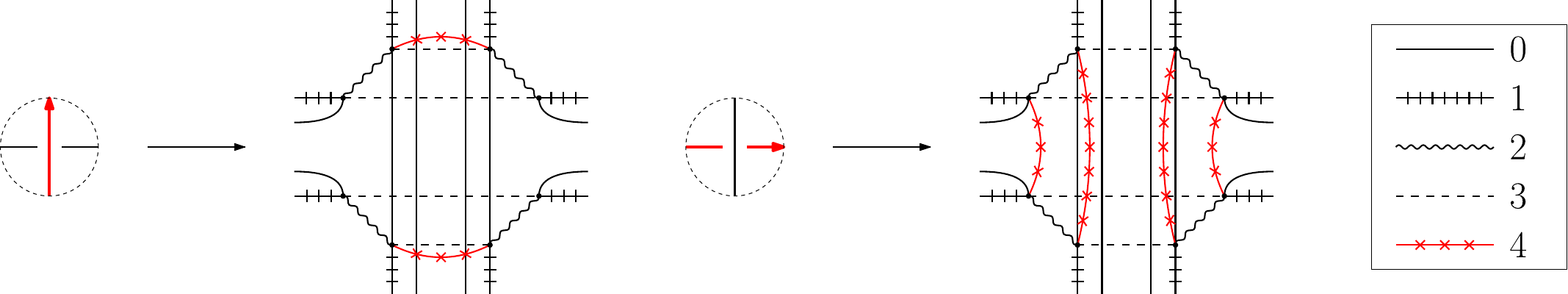}
\caption{Criteria for adding colour-4 edges associated to elements of $Y$.}
\label{fig:dgt-highlightSubgraphs}
\end{figure}

\begin{figure}[h]
	\centering
	\includegraphics[width=1.1\textwidth]{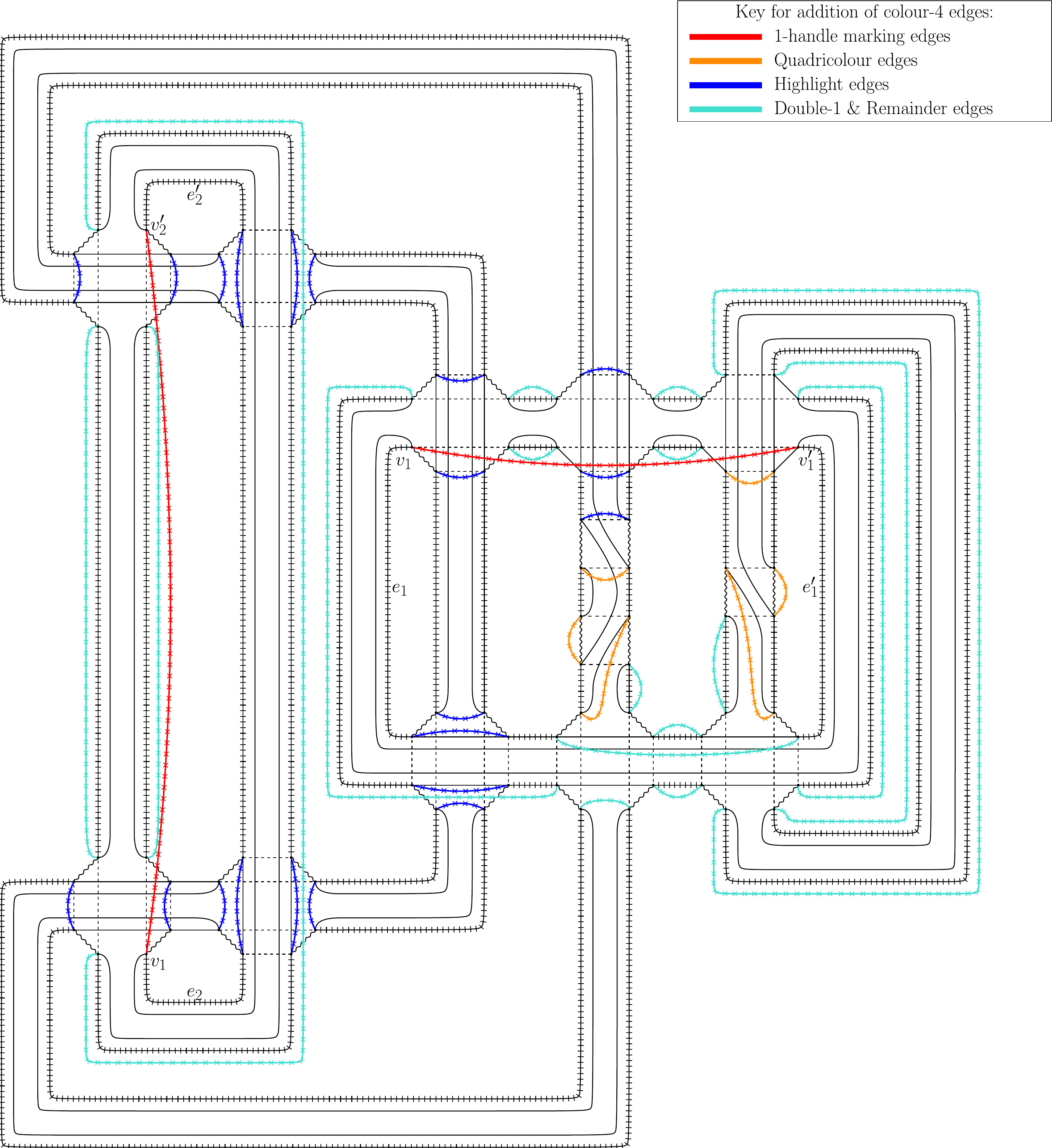}
	\caption{Example of a 5-coloured graph representing the 4-manifold of Figure \ref{fig:kirbyDiagramExample}.}
	\label{fig:dgt-1h-eg-2}
\end{figure}

With the above data in place, the colour-4 edges realising 1- and 2-handle attachments are added according to the following criteria:
\begin{itemize}
	\item For each quadricolour $Q_j$, $j\in\{d+1,\ldots,\ell\}$, add edges according to Figure \ref{fig:quadriSub}.
	\item For each $i\in\{1,\ldots,d\}$, add an edge so as to join $v_i$ and $v_i'$.
	\item For each crossing in $Y_j$, if no colour-4 edges have already been added to the associated subgraph in $\Gamma$, add edges according to Figure \ref{fig:dgt-highlightSubgraphs}. 
	\item Add edges to the remaining vertices, joining those which belong to the same $\{1,4\}$-residue. 
\end{itemize}

Figure \ref{fig:dgt-1h-eg-2} depicts the 5-coloured graph obtained via this procedure when applied to the Kirby diagram of Figure \ref{fig:kirbyDiagramExample}. Note that all edges coloured according to the key are the same colour-4, we have only used different colours to illustrate the separate steps of the procedure.

\clearpage

\section{Up-Down Simplification}
\label{sec:uds}
A drawback of the algorithm used by DGT is that the size of the resulting triangulations are quite large. From Figure \ref{fig:dgt-graphs1}, we see for each true crossing $x$ of the diagram, and each framing curl $c$, the number of pentachora (4-simplices) in the resulting triangulation is $8x+4c$. 

To illustrate this number is not optimal, consider the following simple example. 
The complex projective plane $\mathbb{C}P^2$, admits a handle decomposition consisting of a single 0-, 2- and 4-handle, and can be realised by the Kirby diagram consisting of a single unknot with $+1$ framing. Since DGT uses PD codes to encode link diagram, we have to place a single trivial curl in the unknot to first obtain a diagram which can be encoded by a PD code. We must then add an additional pair of curls to gurantee the existence of a quadricolour. This results in a triangulation of $\mathbb{C}P^2$ with $4\cdot 3=12$ pentachora. However it is not difficult to verify that a minimal triangulation of $\mathbb{C}P^2$ consists of just four pentachora, e.g.\ either by checking the 4-dimensional closed census \cite{regina} or by manually performing Pachner moves (a.k.a bistellar flips, see Section \ref{sec:uds-details}) on the triangulation directly.

\begin{figure}[h]
	\centering
	\includegraphics[width=0.7\textwidth]{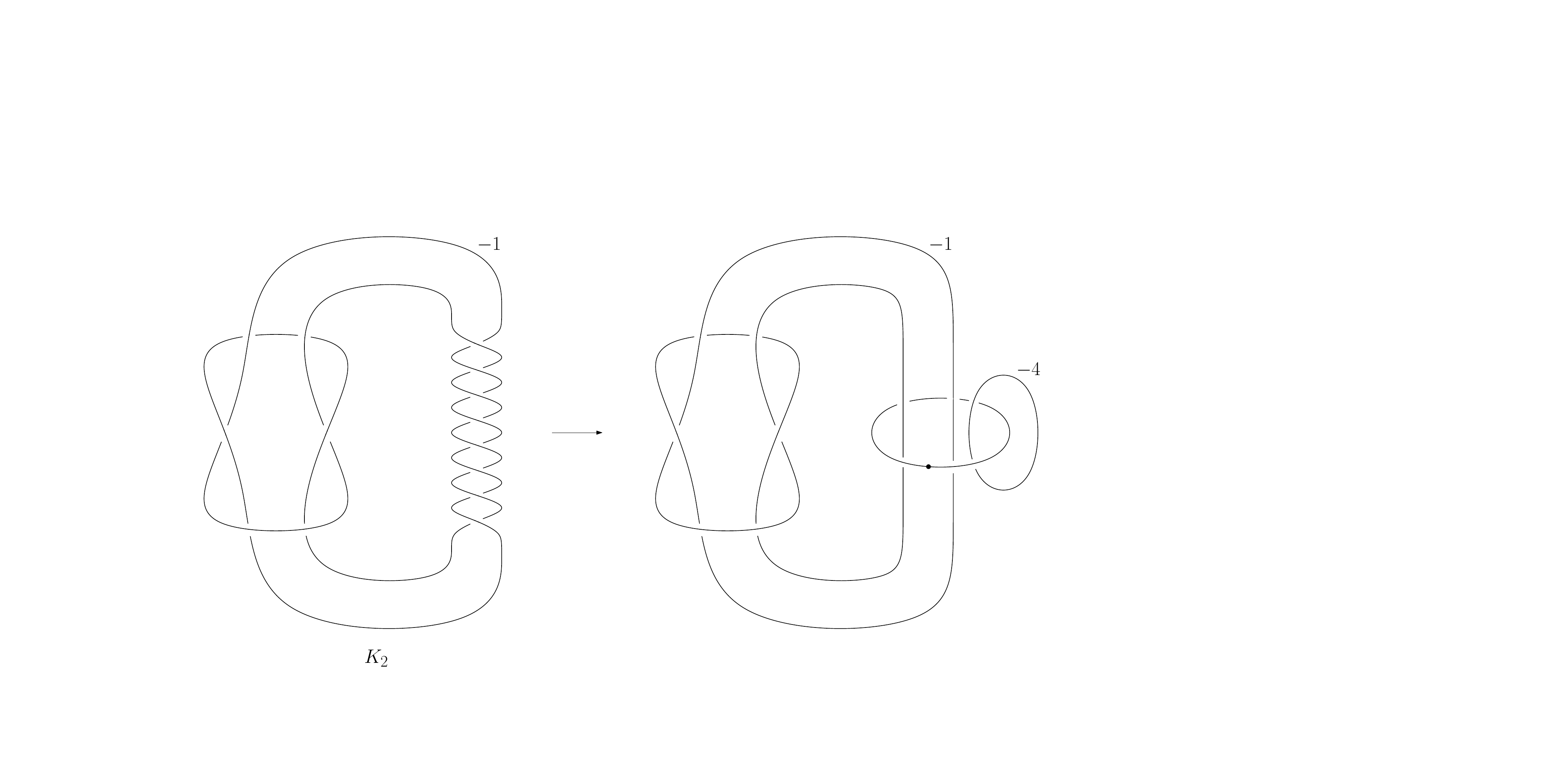}
	\caption{Introducing a cancelling 1/2-handle pair to reduce the complexity of a diagram.}
	\label{fig:A2-12cancelling}
\end{figure}

In some cases introducing a cancelling 1/2-handle pair (see Appendix \ref{app:handleCalculus}) into a diagram can be used to reduce the complexity of the diagram, which in turn reduces the size of the triangulation produced by DGT. For example, consider the diagram on the left in Figure \ref{fig:A2-12cancelling}. By introducing a cancelling 1/2-pair, we can remove the four twists, as seen on the right.

The left diagram had 14 true crossings, with writhe $w(K_2)=+8$ (and so, needing nine additional curls to encode the $-1$ framing). Consequently, the number of pentachora in the triangulation built by DGT is $8\cdot 14+4\cdot 9=148$. If we compare this to the diagram on the right, we see that the cancelling pair has reduced the total number of true crossings in the diagram by two, and in particular, has lowered the writhe of $K_2$ to zero, meaning $K_2$ now only needs one additional curl to encode the framing. Together with the four additional curls to encode the framing of the cancelling unknot, this means that the total number of additional curls needed has also been lowered, from the original nine, down to five. Hence the total number of pentachora in the resulting triangulation is $116$, a saving of 32 pentachora.

Of course, this 1/2 swapping ``trick'' cannot always be performed, and is not guaranteed to lower the complexity when it can be done (for example, compare Figure \ref{fig:ep-gen} to Figure \ref{fig:gen-1hver} in which the 1/2 swap \textit{increases} the complexity of the left link). Moreover even when it does lower the complexity of the link, the triangulations produced by DGT are still fundamentally constrained by the $8x+4c$ ``barrier''. To this end, we have developed a new heuristic called \textit{Up-Down Simplification} (UDS) for reducing the size of 4-manifold triangulations.

\subsection{The Heuristic}
\label{sec:uds-details}
Using the Dehn-Sommerville relations and the equation for the Euler characteristic, it can be shown that the number of pentachora in a triangulation of a simply-connected 4-manifold $M$ is 
\begin{equation}
	p=6(\chi-v)+2e,
	\label{eq:numPent-sc}
\end{equation}
where $\chi$ is the Euler characteristic of $M$, $v$ is the number of vertices, and $e$ is the number of edges. Since $\chi$ is a fixed constant, $p$ is determined by $v$ and $e$. As such, let us define the ``core'' $f$-vector of a 4-manifold triangulation to be $\bar{f}=(v,e\,\vert\,p)$ (recall the standard $f$-vector is the vector whose $i$th entry is the number of $i$-dimensional simplices in a complex).

\begin{figure}[h]
	\centering
	\includegraphics[width=0.75\textwidth]{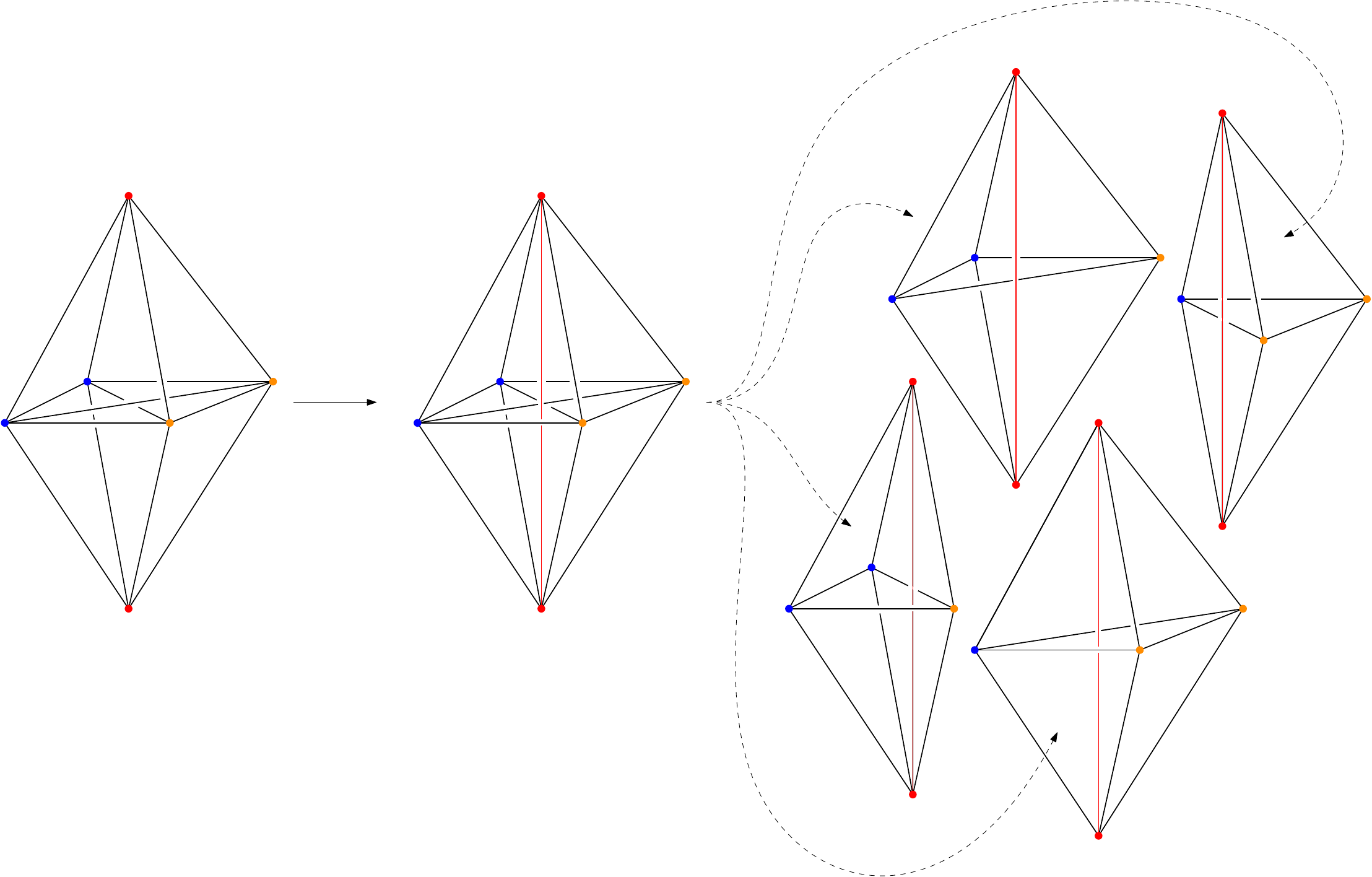}
	\caption{A 2-4 move divides two adjacent pentachora into four.}
	\label{fig:pachner24}
\end{figure}

We modify our triangulations using local moves which change the triangulation but not the underlying topology or PL type, for example Pachner moves (a.k.a bistellar flips) \cite{Pachner}. Informally, an $i$-$j$ Pachner move can be thought of as taking $i$ pentachora and replacing them with $j$ pentachora in a specified manner. An example of a 2-4 move is shown in Figure \ref{fig:pachner24}. 

Moves of particular importance to us are (i) 3-3 and 4-4 moves which do not effect the $f$-vector, and (ii) 2-0 edge/triangle moves, both of which change $\bar{f}$ by $(v,e\,\vert\,p)\mapsto(v,e-1\,\vert\,p-2)$.
%\begin{itemize}
%	\item 3-3 and 4-4 moves which do not effect the $f$-vector, and
%	\item 2-0 edge/triangle moves, both of which change $\bar{f}$ by $(v,e,p)\mapsto(v,e-1,p-2)$.
%\end{itemize}

In three-dimensions, a defining characteristic of minimal triangulations is having as few vertices as possible (indeed in typical cases minimal 3-manifold triangulations will have just a single vertex \cite{JacoRubi-0eff3mflds}). Consequently, existing heuristics for 3-manifolds try to minimise the number of vertices. However, from Equation \eqref{eq:numPent-sc}, smaller values of $p$ are achieved not by $v=1$, but by having $v$ and $e$ of the same order of magnitude as $\chi$. In other words, $p$ is minimised by having as many vertices as possible with as few edges as possible.
Consequently, existing 3-manifold simplification heuristics are not optimal for use on 4-manifold triangulations. 

Moreover, the simplification heuristic in \regina for 3-manifolds does not attempt to make the triangulation any bigger than the initial input, this is because in three dimensions local minima are never deep enough to require ``well climbing'' techniques \cite{Bab-HyamFest}. However in the case of 4-manifolds, by increasing the size of the triangulation in a controlled manner, new simplifying moves are often opened up, and importantly, in sufficient number to escape local minima, which in contrast to three dimensions, are often deep enough to get trapped in. 

With the above in mind, the three key ideas underpinning our heuristic are:
\begin{itemize}
	\item do not perform any moves which change the number of vertices, i.e.\ no 1-5/5-1, 2-0 Vertex, or Collapse Edge moves (for details of these moves see Appendix \ref{app:pachnerMoves} or \cite{Bab-HyamFest}),
	\item only descend via 2-0 Edge, and 2-0 Triangle moves, and
	\item use 2-4 moves to enlarge the triangulation to ``open'' new simplification paths (but keep the number of vertices fixed).
\end{itemize}
Note we do not use 4-2 moves to descend since these simply reverse the 2-4 moves performed and hence do not change our position in the space of triangulations. Whereas by using 2-0 moves, we are genuinely ``going down a different path''.

Let \texttt{upMax} and \texttt{sideMax} be some maximum number of 2-4 and 3-3/4-4 moves to perform. 

The core algorithm of the heuristic runs as follows.
\begin{enumerate}
	\item For $1\leq i\leq$\ \texttt{upMax} perform $i$ random, 2-4 moves.
	\item If no 2-0 Edge or 2-0 Triangles moves are available, perform $j$ random 3-3 moves, or if no 3-3 moves are available, perform $j\leq$\ \texttt{sideMax} random 4-4 moves (and vice-versa).
	\item Perform as many 2-0 Edge and 2-0 Triangle moves as possible.
	\item Repeat from Step 1 as many times as desired.
\end{enumerate}

We note that in practice, once UDS appears to hit a deep local minimum, performing moves which \textit{do} reduce the number of vertices (e.g.\ ``collapse edge'' or 2-0 vertex moves) and then applying UDS to the ``new'', vertex-reduced triangulation can typically be used to some success in further simplifying the triangulation (i.e.\ UDS should be used in conjunction with judicious application of other simplifying moves to get the smallest triangulation possible). Despite the simple ideas behind UDS, it has proven to be more effective than any prior heuristic for 4-manifold triangulations, as demonstrated in the next section.

\section{Experimental Results}
\label{sec:results}
Data pertaining to the triangulations presented in this section are available in Appendix \ref{app:isosigs}.
\subsection{The $K3$ Surface}
\label{sec:small-k3}

The $K3$ surface is one of four ``fundamental'' simply-connected 4-manifolds, alongside $\mathbb{C}P^2$ and $S^2\times S^2$ \cite{MR1707327,Scorpan}. In contrast to $\mathbb{C}P^2$ and $S^2\times S^2$ however which have small Euler characteristics (3 and 4 respectively), $K3$ is comparatively large with $\chi=24$.

Two triangulations of the $K3$ surface are due to Casella and K\"uhnel \cite{K316}, and Spreer and K\"uhnel \cite{K317}. Both triangulations are simplicial complexes, with core $\bar{f}$-vectors $(16,120\,\vert\,288)$ and $(17,135\,\vert\,312)$ respectively (as such, the triangulations are commonly referred to as $K3_{16}$ and $K3_{17}$). Whilst $K3_{17}$ is known to have a standard PL type, the PL type of $K3_{16}$ is unknown. It is conjectured that they are diffeomorphic, and in 2014, Burton and Spreer attempted to show this (unsuccessfully) by adapting methods from three dimensions \cite{BurtonSpreer-K3}. 

An obvious approach to proving the conjecture would be to simplify both triangulations as far as possible, and then repeatedly perform random local modifications, which keep the number of simplices fixed, until both triangulations are combinatorially isomorphic.

As mentioned in Section \ref{sec:uds}, minimal triangulations of 3-manifolds typically have a single vertex, and so Burton and Spreer start by applying this mindset to the triangulations of $K3$, obtaining 1-vertex 1-edge triangulations of $K3_{16}$ and $K3_{17}$ each having $\bar{f}=(1,1\,\vert\,140)$. 

This in itself was a non-trivial task at the time, requiring a combination of techniques including (i) classical techniques which reduce the triangulation as far as possible using local moves, (ii) the ``composite'' moves described earlier, (iii) simulated annealing techniques, and (iv) exhaustive retriangulation in the form of a breadth-first search exploring Pachner moves.
From there, they attempt to connect the two triangulations via the method outlined above, specifically by running a dual-source breadth-first search exploring Pachner moves. At the time, their algorithm detected over 1,738,260 combinatorially distinct 1-vertex 1-edge triangulations of the $K3$ surface. However their algorithm had neither exhausted the list of such triangulations nor connected the two triangulations in question.

Using UDS we obtained triangulations of $K3_{16}$ and $K3_{17}$ with 66- and 60-pentachora respectively. Whilst this is more promising, being considerably smaller triangulations than both the original simplicial complex versions as well as the 1-vertex 1-edge triangulations (as well as being the first triangulations of $K3_{16}$ and $K3_{17}$ with $<100$ pentachora), we have managed to obtain an even smaller triangulation of $K3$ as follows.

Using DGT, we obtained a 2048 pentachora triangulation of the (standard, and hence, diffeomorphic to $K3_{17}$) $K3$ surface from the Kirby diagram shown in Figure 12.17 of \cite{Akb-4mflds}. After applications of UDS, we arrived at a triangulation with just 54 pentachora. This is the smallest known triangulation of $K3$ to date. The dual graph is shown in Figure \ref{fig:K3-54p}, and the isomorphism signature is given in Appendix \ref{app:54pK3}. 

\begin{figure}
	\centering
	\includegraphics[width=0.8\textwidth]{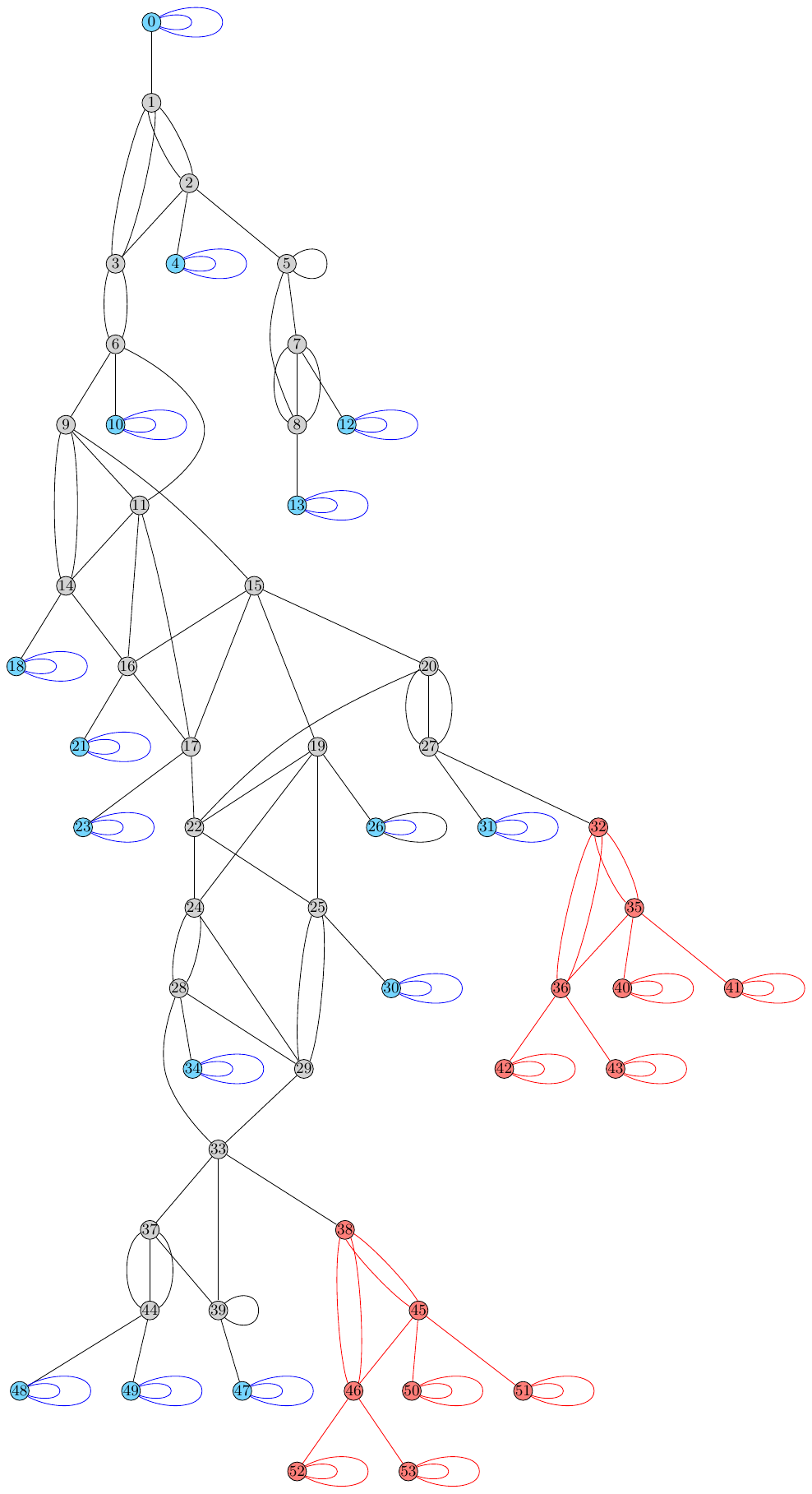}
	\caption{Dual graph of a 54-pentachora triangulation of the $K3$ surface.}
	\label{fig:K3-54p}
\end{figure}

The small size of the triangulation reveals prominent features, which appear to be common to closed simply-connected triangulations. For example, there is a 1-pentachoron triangulation of $B^4$ formed by identifying two pairs of the pentachoron's facets together, which we call a \bemph{twice snapped ball} (TSB). As can be seen in Figure \ref{fig:K3-54p}, this structure appears thoughout the triangulation of $K3$ (nodes highlighted in blue). We believe TSBs play a key role in the handle structure of the manifold (informally, we tend to see one TSB per handle). Similar observations can be made for triangulations of other simply-connected manifolds (e.g.\ $\mathbb{C}P^2$, $S^2\times S^2$, etc.). Moreover, as can be seen in Figure \ref{fig:K3-54p}, these structures appear to ``coalesce'' into larger groups (cf.\ the nodes highlighted in red). It transpires that the 7-pentachora structures indicated all realise the same topological structure, specifically a linear \bemph{plumbing} of three 2-spheres. Plumbings are an important construction in 4-manifold theory, formed by iteratively gluing together disk bundles over surfaces. Many closed simply-connected 4-manifolds contain plumbed manifolds as constituent pieces. 

These observations suggest it may be possible to construct triangulations of simply-connected 4-manifolds in a manner which directly reflects the topological structure of the manifold. In addition, plumbings also play a role in producing exotic structures (via a procedure called \bemph{rational blowdown}, cf.\ \cite{fintushelSternRationalBlowdown,parkRationalBlowdown}), and so being able to recognise these structures is of independent interest in this respect.

Whilst this is the smallest triangulation of $K3$ obtained so far, we believe that an even smaller triangulation is possible:

\begin{conjecture}
	\label{conj:lower-bound-conj}
Let $M$ be a closed, simply-connected, smooth 4-manifold. If $\mathcal{T}(M)$ is a triangulation of $M$, with $p$ pentachora, then $2\chi-2\leq p$.
\end{conjecture}
Evidence for Conjecture \ref{conj:lower-bound-conj} stems from observations of the closed census of up to six pentachora as well as examples constructed via DGT and UDS. This conjecture has also recently been independently investigated by Spreer and Tobin \cite{SpreerTobin2024}. 

If Conjecture \ref{conj:lower-bound-conj} is true in general, then our triangulation of $K3$ is only eight pentachora (or four edge reductions) away from being minimal. Despite minimality seemingly within arm's-reach, we have not been successful in simplifying our triangulation any further, suggesting that this triangulation sits within a very deep local minimum.

\clearpage

\subsection{Exotic 4-Manifolds}
\label{sec:exotic-triangulations}
\begin{figure}[h]
	\centering
	\includegraphics[width=0.7\textwidth]{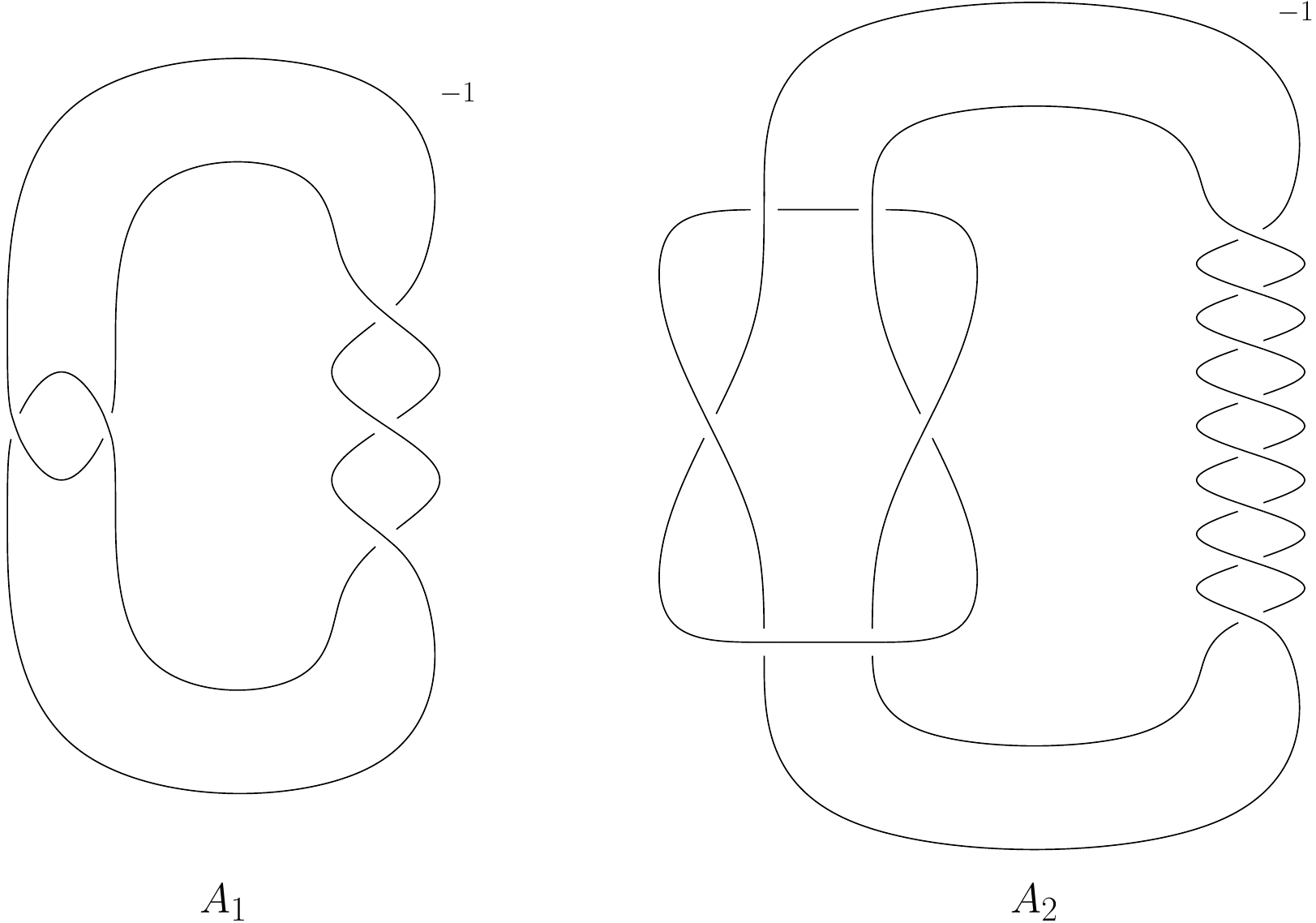}
	\caption{An exotic pair, due to Akbulut \cite{ep-Akb}.}
	\label{fig:ep-akb}
\end{figure}

\begin{figure}[h]
	\centering
	\includegraphics[width=0.9\textwidth]{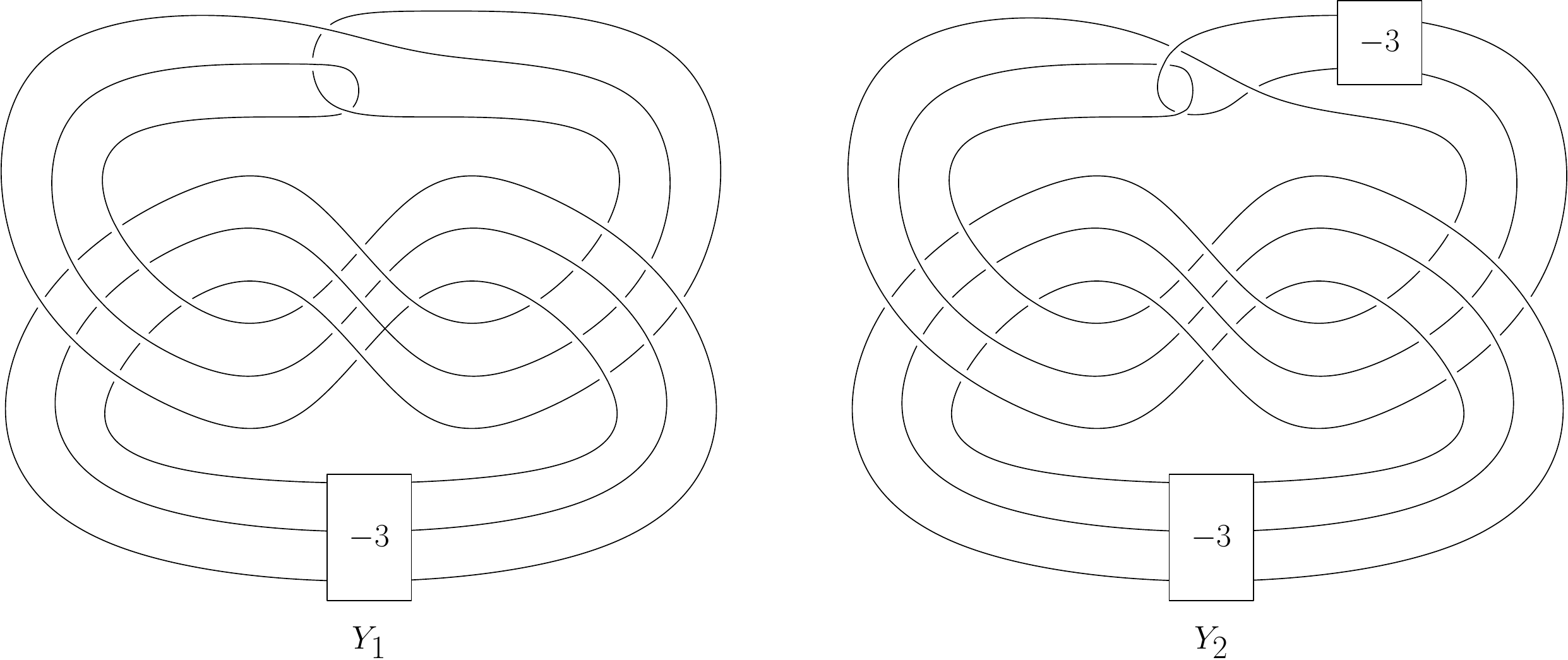}
	\caption{An exotic pair, due to Yasui \cite{ep-Ysi}.}
	\label{fig:ep-ysi}
\end{figure}

\begin{figure}[h]
	\centering
	\includegraphics[width=0.85\textwidth]{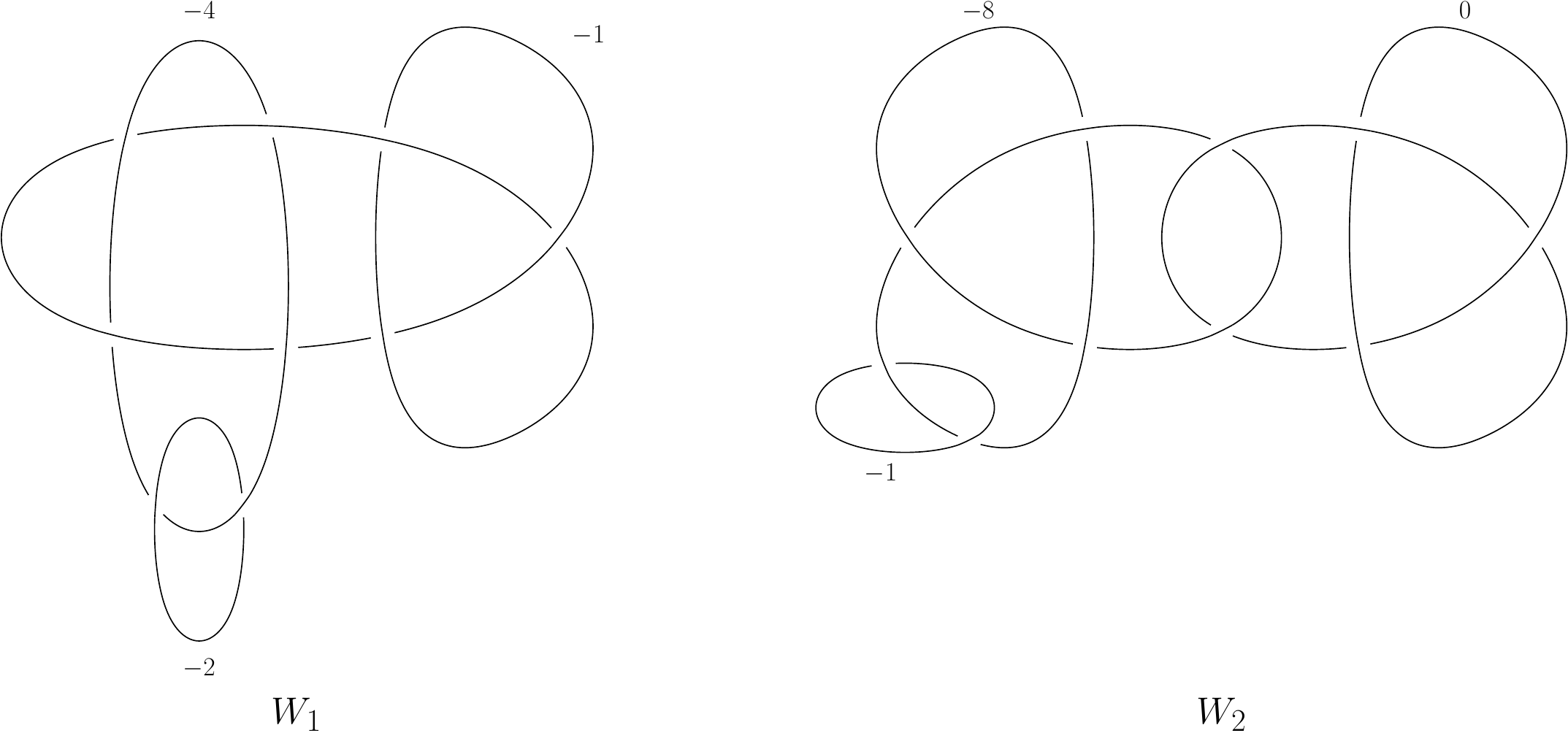}
	\caption{An exotic pair, due to Naoe \cite{ep-Naoe}.}
	\label{fig:ep-naoe}
\end{figure}

\begin{figure}[h]
	\centering
	\includegraphics[width=0.85\textwidth]{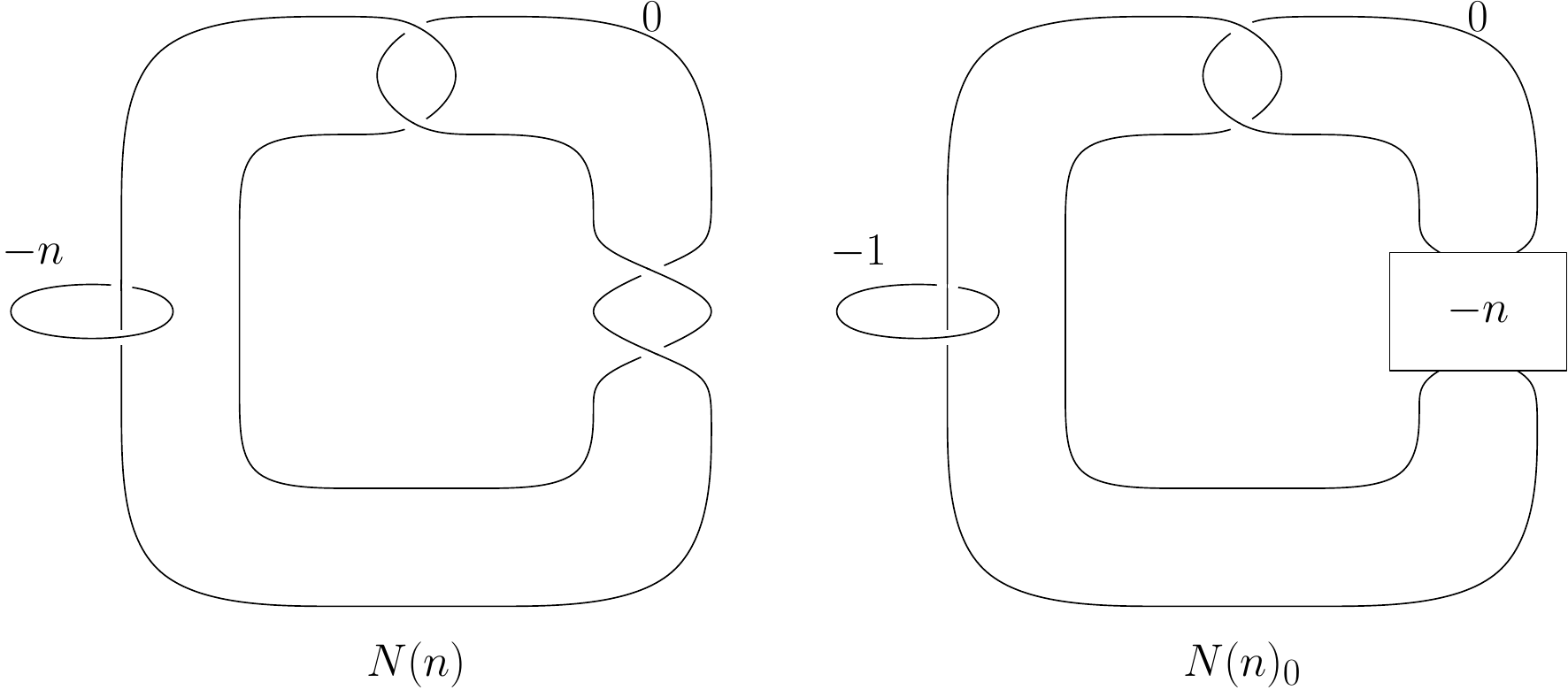}
	\caption{An exotic pair, $n\geq 3$ odd, due to Gompf \cite{ep-Gompf}.}
	\label{fig:ep-gen}
\end{figure}

Using DGT, we produce triangulations of the exotic 4-manifold pairs shown in Figures \ref{fig:ep-akb}--\ref{fig:ep-gen}. These are the first known triangulations of simply-connected orientable exotic 4-manifolds. We note that whilst the 5-colour graphs representing the pair in Figure \ref{fig:ep-naoe} were illustrated in \cite{CasaliCristofori2023Final}, the triangulations presented in this paper constitute the first readily machine-readable presentation of all the pairs in question. Through applications of UDS, we obtain particularly small triangulations, whose $f$-vectors are shown in Table \ref{tab:exotica-fVectors}. The objects with ``cork'' and ``plug'' labels are discussed in Section \ref{sec:corks-and-plugs}.

\begin{table}[h]
	\centering
	\begin{tabular}{c | c || c | c || c | c}
Manifold & $f$-Vector & Pair & $f$-Vector & Pair & $f$-Vector \\
\hline
Akbulut Cork $C_1$ & $(1,2,17,25,10)$ & $(A_1,A_2)$ & $(2,2,17,25,10)$ & $(W_1,W_2)$ &$(4,5,30,40,16)$ \\
Positron Cork $\bar{C}_1$ & $(1,6,34,45,18)$ & $(N(3),N(3)_0)$ &$(3,3,22,30,12)$ & $(M_1,M_2)$ & $(2,9,46,60,24)$ \\
Plug $P_{1,2}$ & $(3,3,12,15,6)$ & $(N(5),N(5)_0)$ & $(3,5,30,40,16)$ & $(Y_1,Y_2)$ & $(2,10,50,65,26)$\\
	\end{tabular}
	\caption{$f$-Vectors of various exotic pairs and related objects.}
	\label{tab:exotica-fVectors}
\end{table}

\begin{remark}
	The exotic pair of Figure \ref{fig:ep-akb} is the smallest and simplest known exotic pair.
\end{remark}

\begin{figure}[h]
	\centering
	\includegraphics[width=\textwidth]{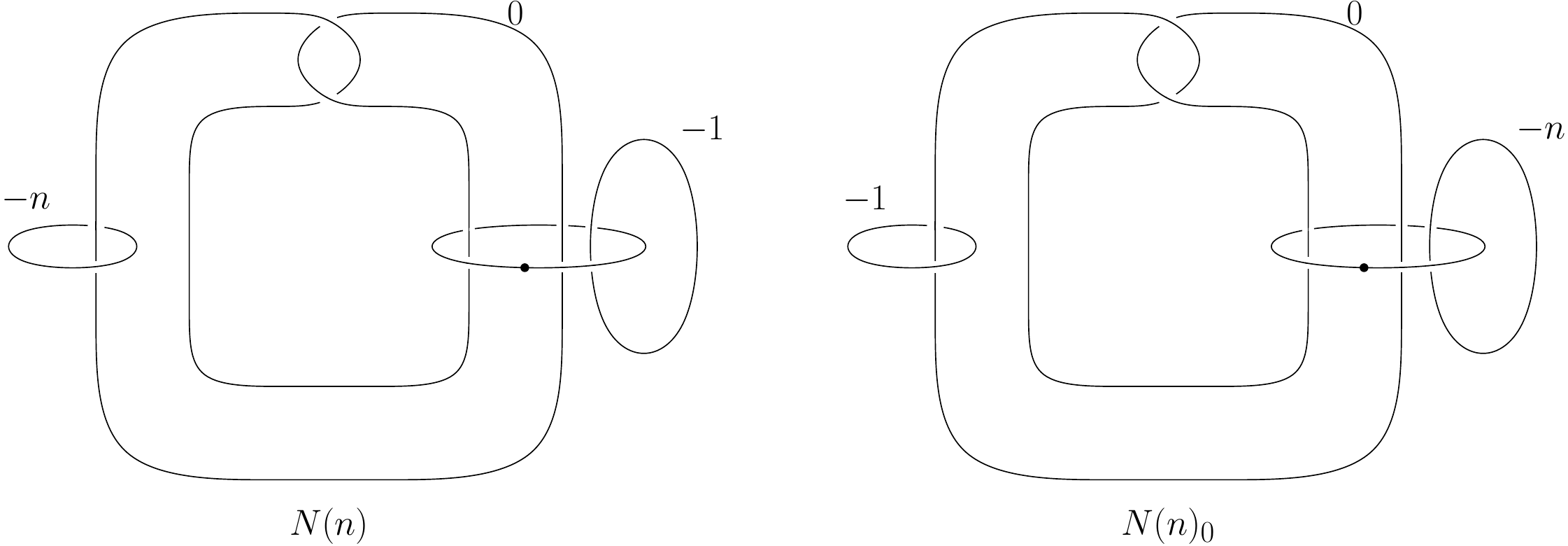}
	\caption{A more symmetric view of the exotic pair from Figure \ref{fig:ep-gen}.}
	\label{fig:gen-1hver}
\end{figure}

Let us consider the pair shown Figure \ref{fig:ep-gen}. By introducing a cancelling 1/2-pair we get the diagrams shown in Figure \ref{fig:gen-1hver}. Observe that both manifolds now have the same underlying link diagram with the only change being which unknot component receives the $-n$ framing. In principle then, it would be possible to automatically generate an infinite family of (exotic) triangulations from this diagram. 

To the best of the author's knowledge, this likely constitutes the simplest such family in the sense that the underlying diagrams of the pair are the same, as well as the links themselves being comparatively simple.

The small size of the triangulations has enabled us in some cases to combinatorially identify a number of distinguishing topological features in some of the exotic pairs. 

For example, both $N(n)_0$ and $W_2$ admit a splitting by a $\overline{\mathbb{C}P^2}$ connect summand. In the Kirby diagrams of the respective manifolds, this is identified by the $(-1)$-framed unknot in each diagram. In particular, this shows that $N(n)_0$ and $W_2$ are not diffeomorphic to $N(n)$ and $W_1$ respectively, since neither of the latter manifolds admit such a splitting \cite{ep-Gompf}.

\begin{figure}[h]
	\centering
	\includegraphics[width=\textwidth]{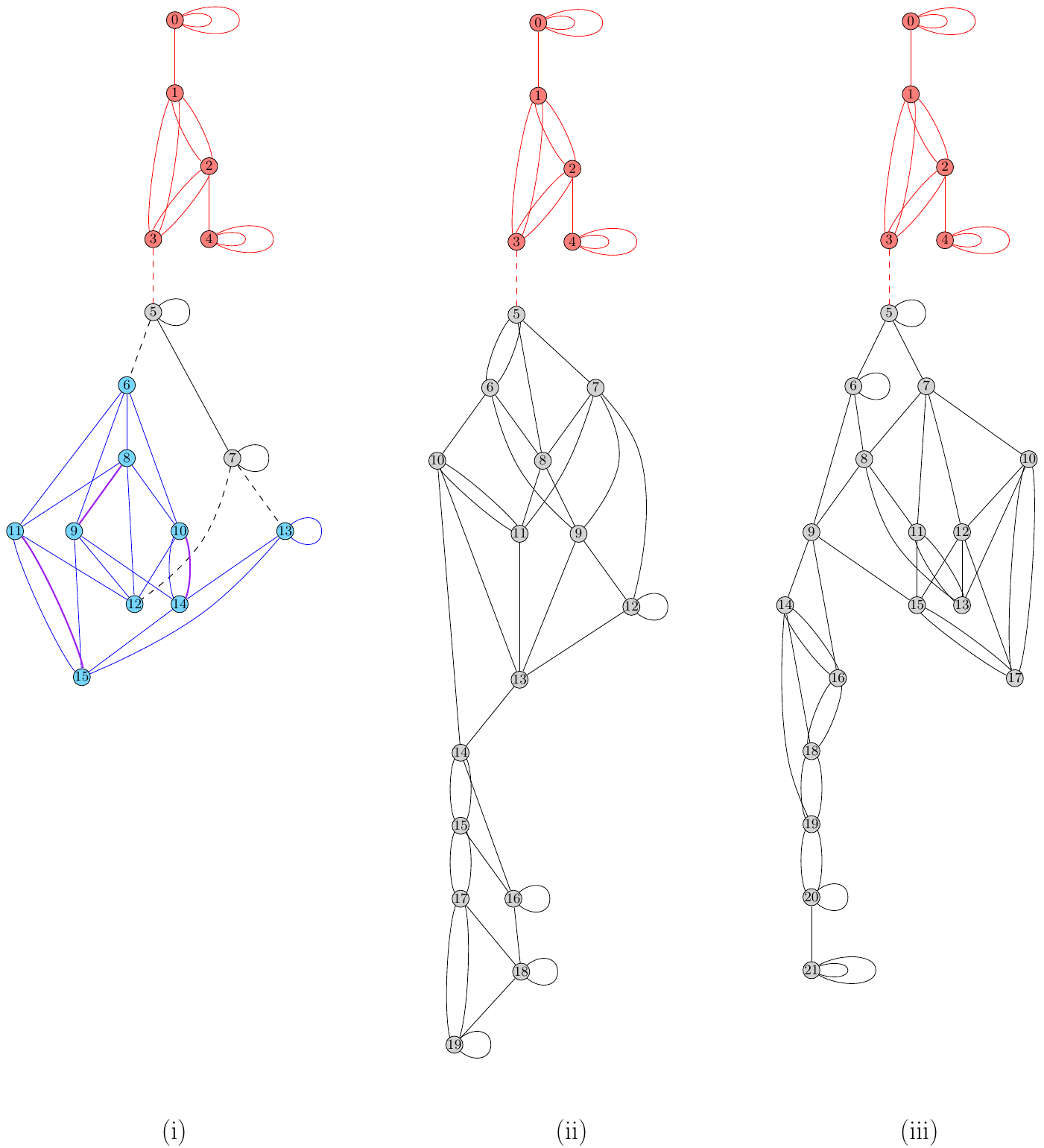}
	\caption{Dual graphs of triangulations of $N(3)_0$ (left), $N(5)_0$ (middle), and $W_2$ (right).}
	\label{fig:N30+N50+W2}
\end{figure}

Using UDS, we obtain a 16 pentachora triangulation of $N(3)_0$, whose dual graph is shown in Figure \ref{fig:N30+N50+W2}(i). If one ``cuts'' the triangulation along the dashed red edge of Figure \ref{fig:N30+N50+W2}(i) (joining pentachora 3 and 5) and fill in the respective $\mathbb{S}^3$ boundary components in each resulting component, then we obtain two new triangulations, say $X$, and a triangulation of $\overline{\mathbb{C}P^2}$ (coming from the subcomplex highlighted in red). In other words, we find the splitting $N(3)_0=X\#\overline{\mathbb{C}P^2}$. In $X$, if one considers the subcomplex derived from only the blue edges of Figure \ref{fig:N30+N50+W2}(i), we get a triangulation of $C(\partial N(3)_0)$ (i.e.\ a \bemph{cone} over the boundary of $N(3)_0$). It would seem then that the additional gluings, indicated by the purple edges, are responsible for realising a cobordism between $C(\partial N(3)_0)$ and $\mathbb{S}^3$. An analogous decomposition is possible for both a 20-pentachora triangulation of $N(5)_0$ (Figure \ref{fig:N30+N50+W2}(ii)) and a 22-pentachora triangulation of $W_2$ (Figure \ref{fig:N30+N50+W2}(iii)), wherein we find the exact same 5-pentachora subcomplex realising the $\overline{\mathbb{C}P^2}$ summand of each manifold. These examples represent a first small step towards being able to combinatorially distinguish smooth structures.

In forthcoming work we will analyse, among others, triangulations of some of these manifolds with real boundary (as opposed to the ideal boundary the present triangulations have) in which the handle structure of the manifold becomes far more apparent --- namely we see a direct ``$D^4\cup 2\text{-handle}$'' structure emerge \cite{RAB_StructurePaper}.

\clearpage

\subsection{Corks and Plugs}
\label{sec:corks-and-plugs}

Given an exotic pair a natural question to ask is: \textit{what causes the change in smooth structure?} In the smooth setting, the following definition and theorem shed some light on this question.

\begin{definition}\label{def:cork}
A cork is a pair $(C,f)$, where $C$ is a compact contractible manifold, and $f:\partial C\to\partial C$ is an involution, which extends to a self-homeomorphism of $C$, but does not extend to a self-diffeomorphism of $C$. We say $(C,f)$ is a cork of a manifold $M$ if $C\subset M$ and removing $C$ from $M$ and regluing via $f$ gives an exotic copy $M'$ of $M$:
\begin{equation}\label{eq:corkDecomp}
	M=N\cup_{\mathrm{id}}C,\quad M'=N\cup_f C,\quad\quad\quad N=M-\mathrm{int}(C),
\end{equation}
\end{definition}

\begin{theorem}[\cite{CorkDecomp-CFHS,CorkDecomp-Mat}]
	\label{thm:cork-thm}
Every exotic simply-connected, closed pair $(M,M')$ satisfy \eqref{eq:corkDecomp}.
\end{theorem}

\begin{figure}[h]
	\centering
	\includegraphics[width=0.8\textwidth]{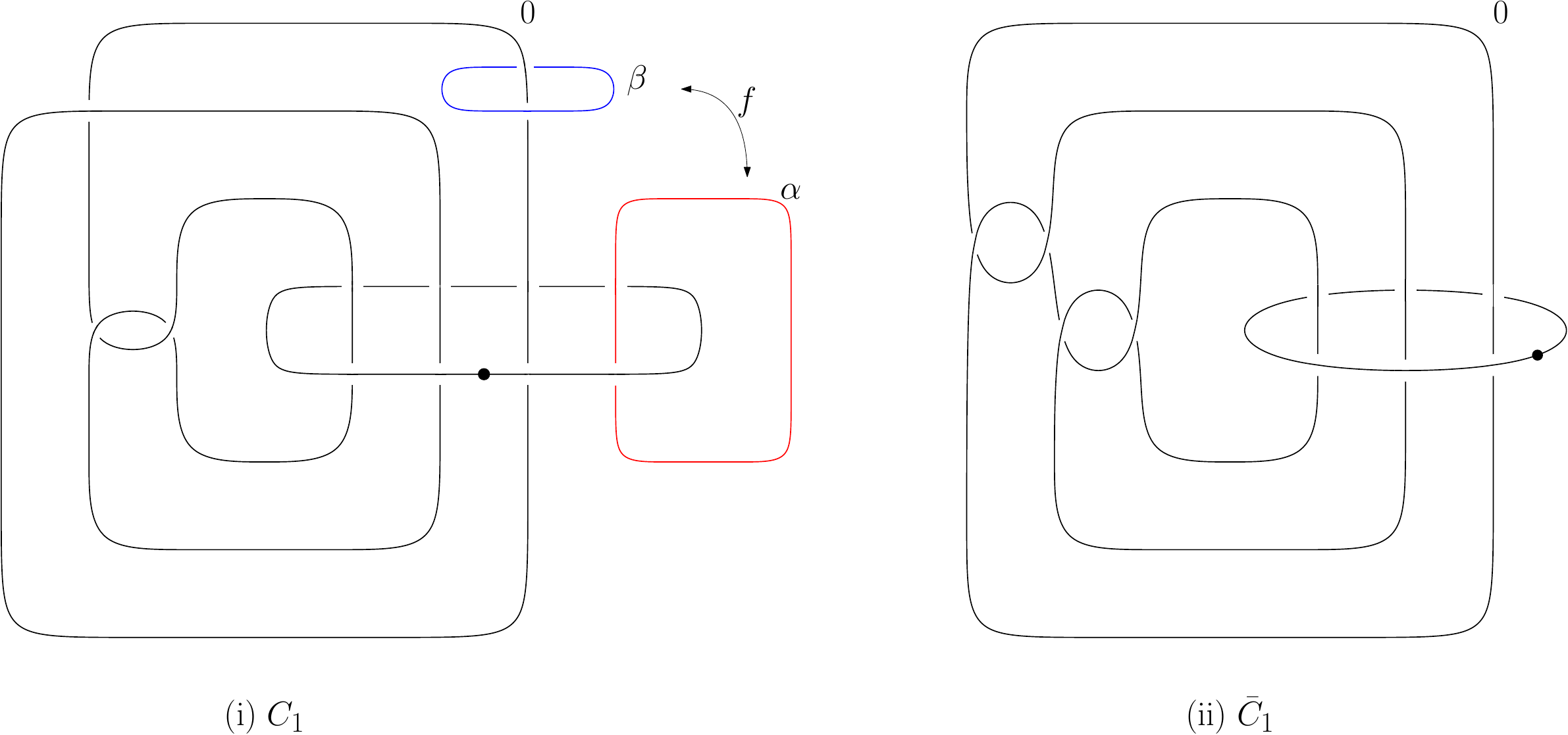}
	\caption{Left: The Akbulut cork $C_1$. Right: The positron cork $\bar{C}_1$.}
	\label{fig:akb-cork}
\end{figure}

Figure \ref{fig:akb-cork}(i) shows the first and simplest known example of a cork, discovered by Akbulut \cite{Akb-Cork}. The map $f:\partial C_1\to\partial C_1$ has the effect of exchanging the loops $\alpha$ and $\beta$.

Many exotic pairs have been constructed from this cork by enlarging it (attaching additional handles). For example, the $(A_1,A_2)$ pair in Figure \ref{fig:ep-akb} is obtained by attaching a ($-1$)-framed unknot along $\alpha$ or $\beta$. See, for example, Figure 11.14 of \cite{MR1707327} for the relevant calculus. Similarly, the $(M_1,M_2)$ pair in Table \ref{tab:exotica-fVectors} are obtained by attaching a $0$-framed trefoil along $\alpha$ (respectively $\beta$)) (see Figure 10.18 of \cite{Akb-4mflds}).

As another example, let $E(1)=\mathbb{C}P^2\#_9\overline{\mathbb{C}P^2}$. The Dolgachev surface is an exotic copy of $E(1)$ \cite{Donaldson-Dolgachev}, and Akbulut showed that the smooth structure can be changed by the so-called ``positron'' cork $\bar{C}_1$ \cite{Akb-Dolgachev}, shown in Figure \ref{fig:akb-cork}(ii).

Once again, using DGT and UDS, we obtain triangulations of the corks $C_1$ and $\bar{C}_1$. This is the first time explicit triangulations of these objects have been obtained.

\begin{figure}[h]
	\centering
	\includegraphics[width=0.45\textwidth]{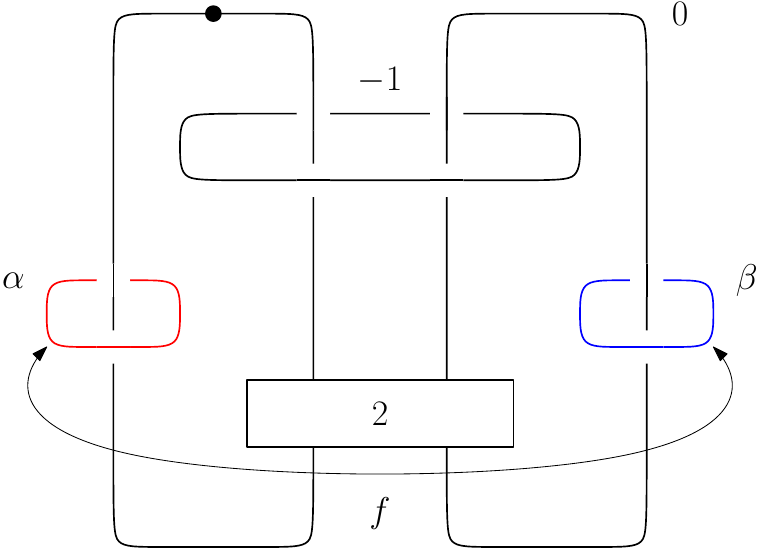}
	\caption{The plug $P_{1,2}$.}
	\label{fig:plug}
\end{figure}

\bemph{Plugs} are similar to corks in that they can sometimes change the diffeomorphism type of a manifold by removing and regluing via an involution on their boundary. Plugs can naturally arise when performing the previously mentioned rational blowdown operation on 4-manifolds \cite{AkbYsi-Plugs}, and so also constitute important objects to have triangulations of, and to study the combinatorics of. Similarly to corks, exotic pairs have been obtained by enlarging plugs, for example the $(W_1,W_2)$ pair of Figure \ref{fig:ep-naoe}, differ by a so-called $P_{1,2}$ plug shown in Figure \ref{fig:plug} (see also Figure 13 of \cite{ep-Naoe} for the relevant handle calculus for $(W_1,W_2)$).

Since corks and plugs are objects that can be responsible for ``exoticity'', having triangulations of these objects will hopefully provide a basis for understanding how triangulations of exotic pairs differ combinatorially. In addition, locating a cork within a 4-manifold is often a highly non-trivial task, and so an algorithm which could identify a combinatorial cork within a triangulation would be of enormous benefit.

\subsection{Work Towards Closed Examples}
\label{sec:ongoing}
Now that we have an effective means of triangulating 4-manifolds, one might hope we can start obtaining triangulations of \bemph{closed} exotic 4-manifolds. Unfortunately, Kirby diagrams of closed exotic manifolds are currently few in number, but worse still is that all of the currently known diagrams of closed examples have their handles drawn in ways incompatible with the algorithm of DGT (recall Figure \ref{fig:dgt-diagramHandleReqs}). In principle, it would be possible to isotope these 1-handles into standard position, however attempts to do this have been unsuccessful due to the high complexity of the diagrams (compare for example Figure 21 of \cite{AkbYsi-Plugs} against the required conditions illustrated in Figure \ref{fig:dgt-diagramHandleReqs}).

On the other hand, there is some hope we might construct examples ``manually'' by decomposing 4-manifolds into pieces which are easier to triangulate, and then gluing them back together appropriately. One possible example of this is the following procedure introduced by Fintushel and Stern \cite{FintushelStern}. Let $K\subset S^3$ be a knot and $\nu(K)$ its tubular neighbourhood. Suppose $X$ contains an embedded torus $T^2$ with trivial normal bundle (i.e.\ $T^2\times D^2\subset X$). Since $S^3-\nu(K)$ is topologically a solid torus $S^1\times B^2$ we have
\[
\partial(T^2\times D^2)=T^2\times\partial B^2=T^2\times S^1=S^1\times\partial B^2\times S^1=\partial(S^3-\nu(K)\times S^1).
\] 
Let $X_K=(X-T\times D^2)\cup_\varphi([S^3-\nu(K)]\times S^1)$, where $\varphi$ is any map which preserves the homology class $[\pt\times\partial D^2]$ of the torus.
If $X$ and $X-T^2\times D^2$ are both simply-connected, then $X_K$ is homeomorphic to $X$, but the diffeomorphism type is determined by $\varphi(\pt\times S^1)$. 

\begin{figure}[h]
	\centering
	\includegraphics[width=0.8\textwidth]{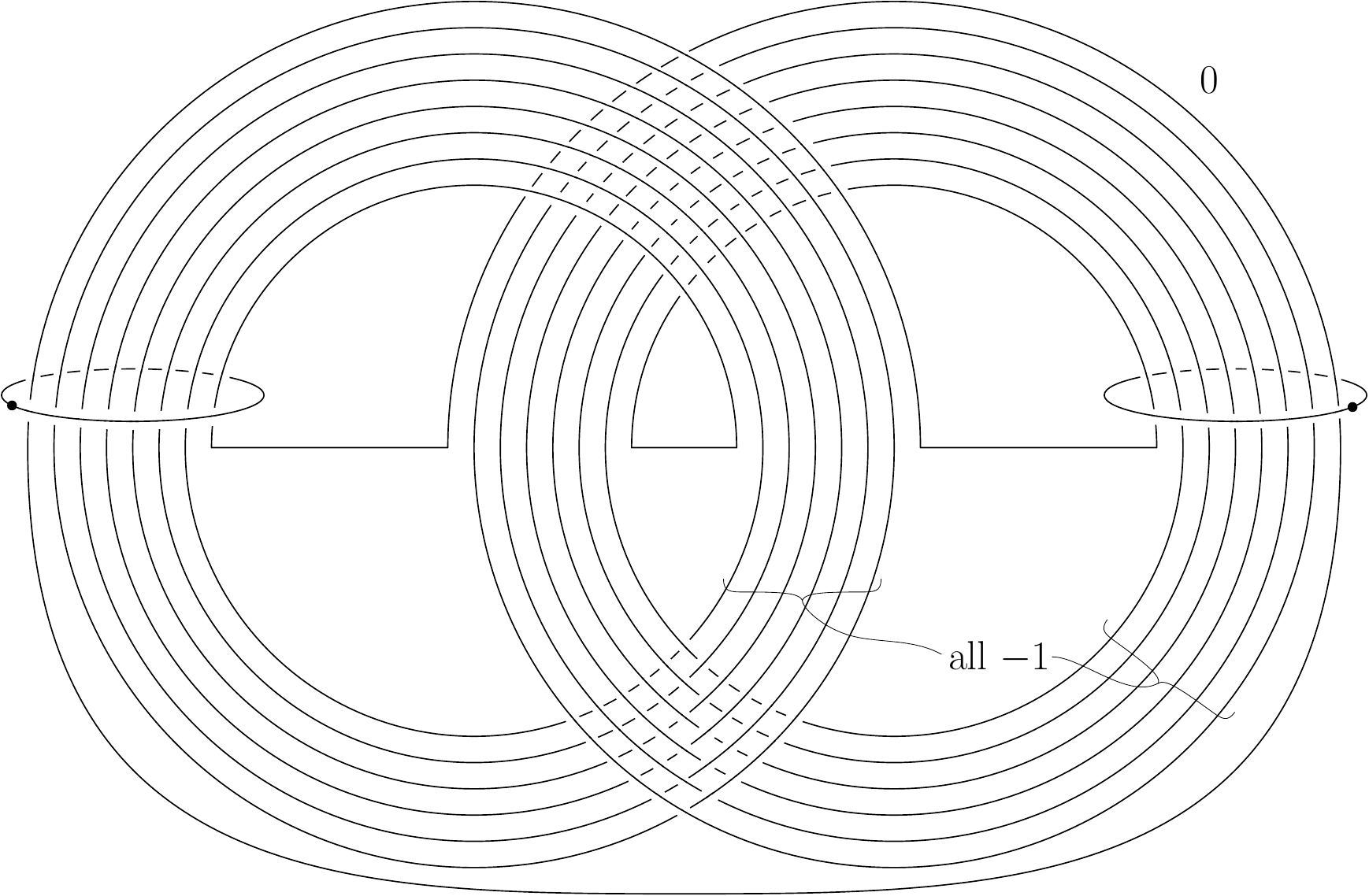}
	\caption{$E(1)-\nu(T^2)$}
	\label{fig:e1-nT2}
\end{figure}

Recall that $E(1)\cong\mathbb{C}P^2\#_9\overline{\mathbb{C}P^2}$. Take $K$ to be the right-handed trefoil knot. It is known that $E(1)_K$ is diffeomorphic to the Dolgachev surface mentioned in Section \ref{sec:corks-and-plugs} (and hence, not diffeomorphic to $E(1)$) \cite{park2001noncomplex}. Using a combination of DGT, UDS, and \regina, we obtained a 40-pentachora triangulation of $E(1)-\nu(T^2)$ from Figure \ref{fig:e1-nT2}. 

Using Regina, we then built a triangulation of $(S^3-\nu(K))\times S^1$ (since \regina already has the capability to build knot complements and $S^1$ bundles of arbitrary 3-manifolds, we did not need to build $(S^3-\nu(K))\times S^1$ from a Kirby diagram). Again using a combination of UDS and \regina we obtained a 16-pentachora triangulation of this piece. We used UDS in conjunction with boundary specific moves to simplify the two triangulations as much as possible, and to the point where the boundaries of both triangulations were the minimal six-tetrahedra $T^3$. The final manifold is then obtained by gluing the boundaries together. There is only a finite number of ways of identifying the boundary triangulations together, and we can perform all such gluings to obtain our candidate exotica.

Whilst the condition on $\varphi$ is not particularly restrictive, it remains to verify that a given gluing of the boundaries does satisfy the necessary conditions to ensure that the resulting space is exotic --- this would entail locating a representative curve for $[\pt\times S^1]$ and ``tracking its journey'' through the construction, which is feasible but which requires a considerable amount of combinatorial ``bookkeeping'' which has yet to be done. However, it seems reasonable that at least one of these gluings should meet the condition. If true, such a triangulation would constitute the first known example of a closed, simply-connected, orientable, exotic triangulation. Additionally, if it transpires that $K3_{16}$ is in fact not diffeomorphic to $K3_{17}$ (Section \ref{sec:small-k3}), then this would in fact constitute the first example of an \emph{irreducible} (i.e.\ not a smooth connect-sum of other manifolds) closed, simply-connected, orientable, exotic triangulation.

\bibliographystyle{plainurl}
\bibliography{references}

\clearpage

\appendix
\section{1-Handle Notation}
\label{app:1handleNotation}
The attaching map for a 1-handle is $S^0\times D^3\to S^3$. So we can visualise a 1-handle as a pair of disjoint 3-balls. Observe that by treating the 0-handle as another copy of $D^1\times D^3$ and identifying a copy of $S^0\times D^3\subset\mathbb{S}^3$, we can see that the attaching map (now of the form $S^0\times D^3\to S^0\times D^3\subset\mathbb{S}^3$), has the effect of identifying the two $D^1$s together along their boundary $S^0$s, giving $S^1$. Hence, the result of attaching the 1-handle is $S^1\times D^3$.

\begin{figure}[h]
	\centering
	\includegraphics[width=\textwidth]{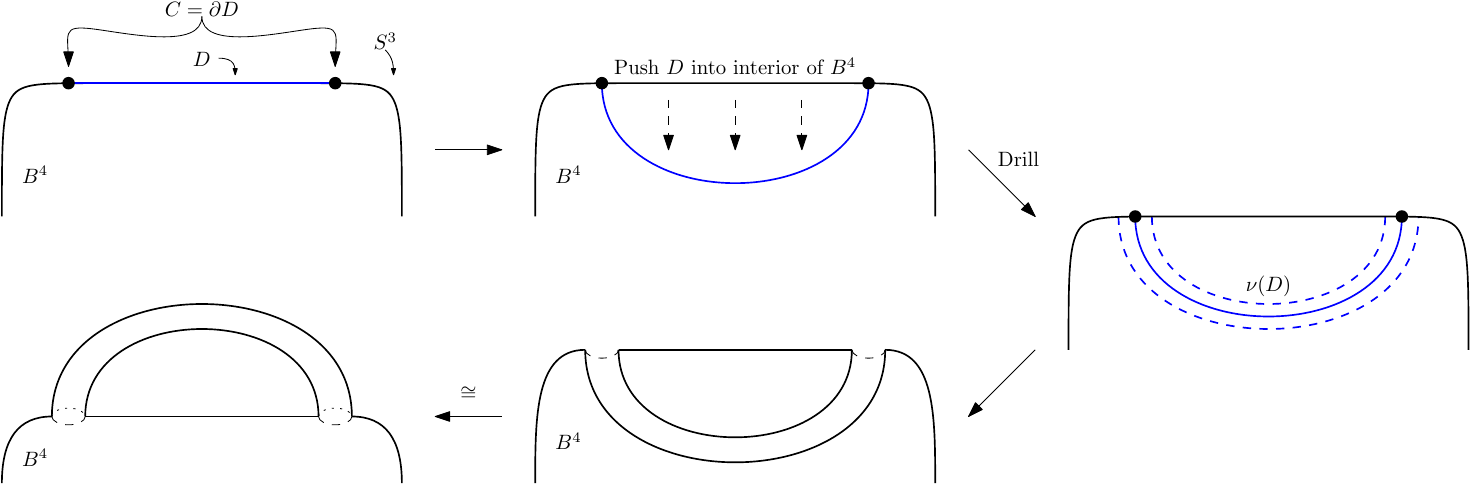}
	\caption{Adding a bridge is the same as removing a tunnel.}
	\label{fig:1h-drilling}
\end{figure}

Certain technical issues can arise however when using the ``disjoint ball'' method of drawing 1-handles, and for this reason we more often use the ``dotted circle'' notation for a 1-handle, which we will now explain. To aide in understanding the following construction, imagine oneself as a topological engineer and chant the mantra ``building a bridge is the same as drilling a tunnel''. First, consider an unknotted circle $C$ in $S^3$. This circle bounds a disk $D$. Push $D$ into the interior of the 4-ball, and then remove a tubular neighbourhood of it. The remainder is $S^1\times D^3$ --- i.e.\ the result of attaching a 1-handle (Figure \ref{fig:1h-drilling}). We can therefore visualise a 1-handle by drawing the boundary of this carved disk, and we place a dot on it to distinguish it from a 2-handle. 

\section{Handle Calculus}
\label{app:handleCalculus}

\begin{figure}[h]
	\centering
	\includegraphics[width=0.95\textwidth]{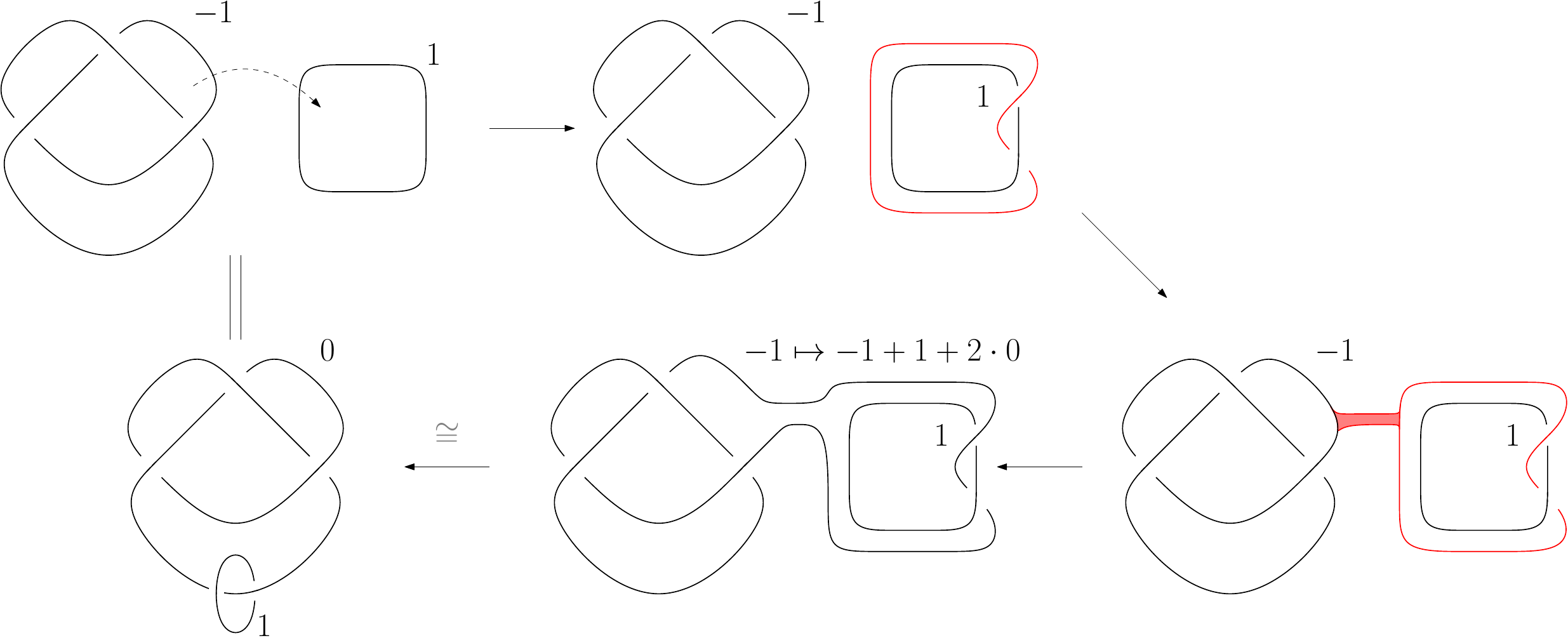}
	\caption{Sliding 2-handles.}
	\label{fig:2handleSlide}
\end{figure}

In \cite{MR0467753}, Kirby introduces a set of moves which can be performed on link traces such that two framed link diagrams represent the same 4-manifold (up to diffeomorphism) if and only if one diagram can be obtained from the other through a sequence of these moves. One of these moves is the action of sliding one 2-handle over another. Diagrammatically this is performed as follows. Given two components $K_i$ and $K_j$ of the link, with framings $\alpha_i$ and $\alpha_j$ respectively, use the framing of $K_j$ to push off a parallel copy of itself, $K'_j$. Replace $K_i$ by $K_i\#_b K'_j$, where $b$ is a band connecting $K_i$ and $K'_j$. The framing of $K_i$ changes according $\alpha_i+\alpha_j\pm 2\lk(K_i,K_j)$, where $\lk(K_i,K_j)$ is the linking number of $K_i$ and $K_j$. This operation is demonstrated in Figure \ref{fig:2handleSlide}.

\begin{figure}[h]
	\centering
	\includegraphics[width=\textwidth]{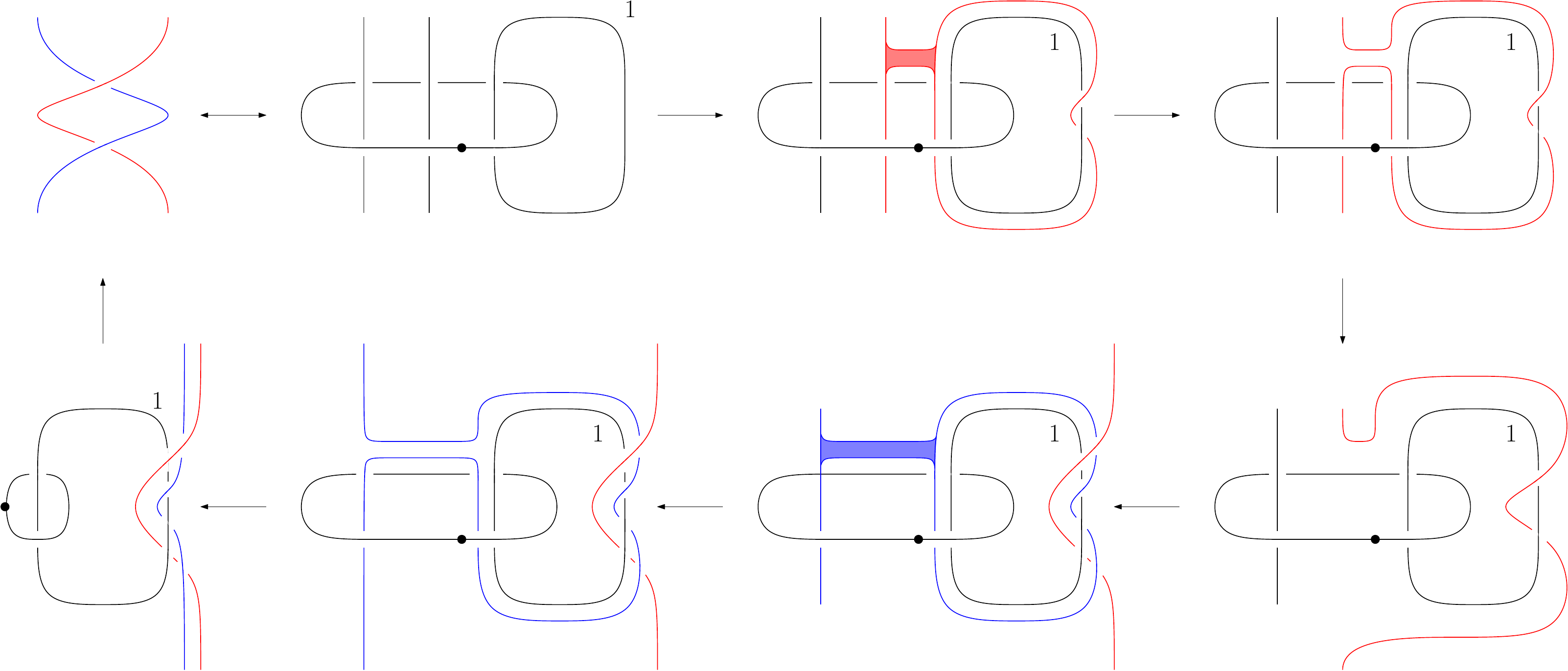}
	\caption{Introducing a cancelling 1/2-pair to remove twists across strands.}
	\label{fig:twist4cancellingPair}
\end{figure}

Another concept which is important is the idea of cancelling 1- and 2-handles. If the attaching circle of a 2-handle intersects the bounding circle of a 1-handle transversally in a single point, then that pair of handles can be erased from the diagram. If there are any other 2-handles passing over the 1-handle, then we first slide those 2-handles over the 2-handle of the cancelling pair. We refer the reader to, for example, \cite{MR1707327}, for full mathematical details. This procedure is illustrated in Figure \ref{fig:twist4cancellingPair}. Whilst one would typically want to cancel any and all extraneous handles in order to obtain as ``efficient'' a handle decomposition as possible, it can sometimes be advantageous to purposefully introduce cancelling pairs and manipulate the diagram in order to show a certain property of the manifold or obtain a diagram in a certain form. For example, Figure \ref{fig:twist4cancellingPair} also shows how twists across strands of 2-handles can be undone by introducing a cancelling 1/2-pair.

\clearpage

\section{Pachner Moves}\label{app:pachnerMoves}
Here we give a very brief description of some four-dimensional local moves for generalised triangulations and how they effect the $f$-vector of the triangulation. For further details and illustrations (albeit for the three-dimensional case), see \cite{Bab-HyamFest}.

\begin{itemize}
	\item 5-1 Move: Replaces five pentachora meeting at a degree 5 vertex with a single pentachoron. Changes the core $\bar{f}$-vector by $(v-1,e-5\svert p-4)$.
	\item 4-2 Move: Replaces four pentachora joined along a degree 4 edge with two pentachora glued together along a tetrahedron. Changes $\bar{f}$ by $(v,e-1\svert p-2)$.
	\item 3-3 Move: Replaces three pentachora joined along a degree 3 triangle with three pentachora joined along a new degree 3 triangle that is aligned in a different direction. This has no effect on the $f$-vector.
	\item 2-4 Move: Replaces two pentachora glued together along a tetrahedron with four pentachora joined along a new degree 4 edge. Changes $\bar{f}$ by $(v,e+1\svert p+2)$.
	\item 1-5 Move: Replaces one pentachoron with five pentachora that meet at a new internal degree 5 vertex. Changes $\bar{f}$ by $(v+1,e+5\svert p+4)$.
	\item 4-4 Move: Replaces four pentachora joined along a degree 4 edge with four pentachora joined along a new degree 4 edge that is aligned in a different direction. This has no effect on the $f$-vector.
	\item 2-0 Move (Triangle): Takes two pentachora joined along a degree 2 triangle and squashes them flat. Changes $\bar{f}$ by $(v,e-1\svert p-2)$.
	\item 2-0 Move (Edge): Takes two pentachora joined along a degree 2 edge and squashes them flat. Changes $\bar{f}$ by $(v,e-1\svert p-2)$.
	\item 2-0 Move (Vertex): Takes two pentachora that meet at a degree 2 vertex and squashes them flat. Changes $\bar{f}$ by $(v-1,e-4\svert p-2)$.
	\item Collapse Edge: Takes an edge between two distinct vertices and collapses it to a point. Any pentachora that contained the edge will be ``flattened away''. Changes $v\mapsto v-1$ and decreases $p$ (variable).   
\end{itemize}

\clearpage

\section{Isomorphism Signatures}
\label{app:isosigs}
To reconstruct the triangulations presented in this paper, download \regina, and produce a new 4-manifold triangulation by selecting type of triangulation ``From isomorphism signature'' and pasting in one of the strings given below. Alternatively, go to the following URL to download a \regina data file containing all of the triangulations presented in this paper: \url{https://github.com/raburke/socg24}.

\subsection{54-Pentachora Triangulation of the $K3$ Surface}
\label{app:54pK3}
\is{2ALALAwLPvQAPQzzwAzPQQMLQvPQwALMvPMPPQQvzPQQQQaacddgeefiiilokkommnnqrqrsswBwBvvyzxxCDDAADHEEFFKJIIKSNSTUNOOPPQQRRUVVWWXXYYZZ0011qbSazbzbQbQaJafaububqbubLbDbMbcaDbPboa2abaNaka+aoaGacanaIaLbaajabaLaqaqbubhb3azb+aoafahbnafa6abaKbIaPbTbmbSb1bSbabQaya6abaGacavawayaNaBbyafaMbcaubga6aba2boagaja2aKa}

\subsection{Exotic Pairs}
\subsubsection{Akbulut Pair $(A_1,A_2)$}
\noindent\bemph{$A_1$:}

\is{kLLAvAQQQccdddhihjigjgjiijaaQb6aQbaarb3aaababa6ababaMbqbPb}

\noindent\bemph{$A_2$:}

\is{kLLAvMQQQcceddhfghihijijjjaaba6aQbaaaa4a4aqbtbqb5a5a5acatb}

\subsubsection{Gompf Exotic Nuclei $(N(n),N(n)_0)$}
\noindent\bemph{$N(3)$:}

\is{mLLALQzLQQQcddeffghfhlkijjjkllaa6a6atb3aearbRbobQb3ajb2aob3aPbfaYaMb}

\noindent\bemph{$N(3)_0$:}

\is{mAvLAAMQLQQaadfgehiihihkkkjlllGacaYaYaaaKaobAaiakbiaya8aMaJbga3aubYa}

\noindent\bemph{$N(3)_0$ (16-pentachora version with ``cut edge'' $\overline{\mathbb{C}P^2}$ summand):}

\is{qALAMAvLLQLQQQcaacdddeefhjklmmmooppmnoppqbvaobtayaSa2ava2a2a2a6atb5aqb2arbrbRbRbPbubpadbQb}

\noindent\bemph{$N(5)$:}

\is{qLAAMAwLMLMAQQcbcbcdefefgkjkjlnmnmnnppopMaJbMalbaaja3bPbWbzaaaIbaabaaaIbIbyalb2bNaaaaababa}

\noindent\bemph{$N(5)_0$:}

\is{qAvLALQMLQQzMQcaafdfihfihilllkjmlnmpooppqbga+aKaGaiaMacaiaWa+aiaiaWaJboaAaYawbgawbvaGaaaGa}

\noindent\bemph{$N(5)_0$ (20-pentachora version with ``cut edge'' $\overline{\mathbb{C}P^2}$ summand):}

\is{uALAMMvMLAQPQLPzQkaacdddeegiijljmnllnmnopqrqststt2aoa7apaTb2bPboaVbJaJaDa7aDaBaWaRaRaJb2aPafadb9aVbgaoaVbJaaaJa}

\subsubsection{Naoe Pair $(W_1,W_2)$}
\noindent\bemph{$W_1$:}

\is{qALLPvLPQLQMQQcaabdffkghhhkomolnponpnpopyaGaGaKaobwbobyagaKaubGaaaibwboawbaaaawbaaobobwb+a}

\noindent\bemph{$W_2$:}

\is{qAvLAAMQQLLMMQcaafegegfhhhjjiillnmnopopp+aoa+aGaaaYaobibMaMayaGagayacaububAaibububabYaKbib}

\noindent\bemph{$W_2$ (22-pentachora version with ``cut edge'' $\overline{\mathbb{C}P^2}$ summand):}

\is{wALAMALvzAvQQQPPQMMcaacdddeefgijlnmrpnnnprqqrrssttuuvv2ava4auayaRaqbva2afaabqbtbwb2aDaWa4a4awbWaVbRaRayayaEaVbRaRaVbRayafa}

\subsubsection{Akbulut-Yasui (Cork+Tref(0)) $(M_1,M_2)$}
\noindent\bemph{$M_1$:}

\is{yLLLLMQQQzvvPMAMLQQQQbdedhgfgiihihjjkqqnsptuqwvuxxwwvvvwxxaaSalbNarbaazaaavavaSaIbSaNa2b3baaaaba3bMaPbaaIbJbaaEbMaMaJbJbpaJbJbEblbqb}

\noindent\bemph{$M_2$:}

\is{yLLALwzQPQLPAAwMMMPQQceccgggkjkiimllnnmmmooqsusutuxvvwwxxxaapbMbpbaaaasbpb5apbQbQbNbNbNbRbRbRb6a6a2bYa7aaaaaqbEbra7aObTaTaSaGaLbrara}

\subsubsection{Yasui Pair $(Y_1,Y_2)$}
\noindent\bemph{$Y_1$:}

\is{ALLAwvwMQLQMAMPAvPPQQQkccdedhmkihhoompppoprqosqqyyxxwyxvvwzyzzzaaFbIbFbIbJbaaFboaqbqbkbkbpazazazazauaaaPbbaaababaObObObObrbrbaababaaarbObbabaPb}

\noindent\bemph{$Y_2$:}

\is{ALLALAvMQAAPAzPQLvQMQQkcedchgglililiikkmooorttstqvyywvvxzzxxzyzaaPb6a5aObrbrbaaMbNb6aNbfafacacaobqbqbqbaaaaeaobeadaaaObObdaeaQbMbaaaapbpbsbpbsb}

\subsection{Corks and Plugs}
\subsubsection{Akbulut Cork}
\is{kLLLMAQMQccefggghihfijjjjiEabaaapapaaaAahbWbhadaTaJbDb3boa}

\subsubsection{Positron Cork}
\is{sLLLAvMwwQLQQQQQceffdillkpnnlorjppmmnrrqqoqraaNauaaacahbaaNaNaqauaua2bZaaacanb0aJbJbua2bLbjaCabaBaPb}

\subsubsection{Plug $P_{1,2}$}
\is{gLLAQQccddddfeff8a+a+aGayayaaaga8aca}
\end{document}